\documentclass[12pt]{article}
\usepackage{graphicx}
\usepackage{amsmath,amsthm,amssymb,enumerate}
\usepackage{euscript,mathrsfs}
\usepackage[left=2cm,right=2cm,top=3.5cm,bottom=3.5cm]{geometry}
\usepackage{color}

\usepackage{soul}

\catcode`\@=11 \@addtoreset{equation}{section}

\catcode`\@=12

\allowdisplaybreaks

\newtheorem{Theorem}{Theorem}[section]
\newtheorem{Proposition}[Theorem]{Proposition}
\newtheorem{Lemma}[Theorem]{Lemma}
\newtheorem{Corollary}[Theorem]{Corollary}

\theoremstyle{definition}
\newtheorem{Definition}[Theorem]{Definition}

\newtheorem{Remark}[Theorem]{Remark}

\newcommand{\bTheorem}[1]{
\begin{Theorem} \label{T#1} }
\newcommand{\eT}{\end{Theorem}}

\newcommand{\bProposition}[1]{
\begin{Proposition} \label{P#1}}
\newcommand{\eP}{\end{Proposition}}

\newcommand{\bLemma}[1]{
\begin{Lemma} \label{L#1} }
\newcommand{\eL}{\end{Lemma}}

\newcommand{\bCorollary}[1]{
\begin{Corollary} \label{C#1} }
\newcommand{\eC}{\end{Corollary}}

\newcommand{\bRemark}[1]{
\begin{Remark} \label{R#1} }
\newcommand{\eR}{\end{Remark}}

\newcommand{\bDefinition}[1]{
\begin{Definition} \label{D#1} }
\newcommand{\eD}{\end{Definition}}

\newcommand{\Del}{\Delta_x}

\newcommand{\vuB}{\vc{u}_b}

\newcommand{\tvm}{\widetilde{\vc{m}}}
\newcommand{\tS}{\widetilde{S}}
\newcommand{\bfphi}{\boldsymbol{\varphi}}

\newcommand{\bFormula}[1]{
\begin{equation} \label{#1}}
\newcommand{\eF}{\end{equation}}

\newcommand{\Ov}[1]{\overline{#1}}

\newcommand{\DC}{C^\infty_c}

\newcommand{\aleq}{\stackrel{<}{\sim}}

\newcommand{\vr}{\varrho}

\newcommand{\tvr}{\tilde \vr}
\newcommand{\tvu}{{\tilde \vu}}
\newcommand{\tvt}{\tilde \vt}

\newcommand{\vt}{\vartheta}
\newcommand{\vu}{\vc{u}}
\newcommand{\vm}{\vc{m}}

\newcommand{\vc}[1]{{\bf #1}}

\newcommand{\Div}{{\rm div}_x}
\newcommand{\Grad}{\nabla_x}

\newcommand{\dx}{\,{\rm d} {x}}

\newcommand{\dt}{\,{\rm d} t }

\newcommand{\intO}[1]{\int_{\Omega} #1 \ \dx}

\newcommand{\D}{{\rm d}}

\newcommand{\ep}{\varepsilon}

\def\softd{{\leavevmode\setbox1=\hbox{d}%
          \hbox to 1.05\wd1{d\kern-0.4ex{\char039}\hss}}}
\definecolor{Cgrey}{rgb}{0.85,0.85,0.85}
\definecolor{Cblue}{rgb}{0.50,0.85,0.85}
\definecolor{Cred}{rgb}{1,0,0}
\definecolor{fancy}{rgb}{0.10,0.85,0.10}

\newcommand\Cbox[2]{%
    \newbox\contentbox%
    \newbox\bkgdbox%
    \setbox\contentbox\hbox to \hsize{%
        \vtop{
            \kern\columnsep
            \hbox to \hsize{%
                \kern\columnsep%
                \advance\hsize by -2\columnsep%
                \setlength{\textwidth}{\hsize}%
                \vbox{
                    \parskip=\baselineskip
                    \parindent=0bp
                    #2
                }%
                \kern\columnsep%
            }%
            \kern\columnsep%
        }%
    }%
    \setbox\bkgdbox\vbox{
        \color{#1}
        \hrule width  \wd\contentbox %
               height \ht\contentbox %
               depth  \dp\contentbox
        \color{black}
    }%
    \wd\bkgdbox=0bp%
    \vbox{\hbox to \hsize{\box\bkgdbox\box\contentbox}}%
    \vskip\baselineskip%
}

\date{}

\begin{document}


\title{Navier--Stokes--Fourier system with general boundary conditions}

\author{Eduard Feireisl
\thanks{The work of E.F. was supported by the
Czech Sciences Foundation (GA\v CR), Grant Agreement
18-12719S.}
\and Anton\' \i n Novotn\' y { 
}
}


\maketitle

\centerline{Charles University, Faculty of Mathematics and Physics, Mathematical Institute}
\centerline{Sokolovsk\'a 83, CZ-186 75 Prague 8, Czech Republic}

\centerline{Institute of Mathematics of the Academy of Sciences of the Czech Republic;}
\centerline{\v Zitn\' a 25, CZ-115 67 Praha 1, Czech Republic}

\centerline{feireisl@math.cas.cz}

\centerline{and}

\centerline{IMATH, EA 2134, Universit\'e de Toulon,}
\centerline{BP 20132, 83957 La Garde, France}
\centerline{novotny@univ-tln.fr}

\begin{abstract}

We consider the Navier--Stokes--Fourier system in a bounded domain $\Omega \subset R^d$, $d=2,3$, with physically realistic
in/out flow boundary conditions. We develop a new concept of weak solutions satisfying a general form of relative energy
inequality. The weak solutions exist globally in time for any finite energy initial data and comply with the weak--strong uniqueness principle.

\end{abstract}


\bigskip


\section{Introduction}
\label{I}

Turbulent phenomena in fluid flows that persist in the long run are usually attributed to the interaction of the
system with the outer world through the physical boundary of the fluid domain $\Omega \subset R^d$. Still the overwhelming
majority of theoretical work {on the existence of (weak) solutions of fluid systems and their} long time behavior  concerns models with homogeneous {or periodic} boundary conditions. In the framework of viscous fluids, the most popular is the \emph{no--slip boundary condition} for the fluid velocity $\vu$,
\[
\vu|_{\partial \Omega} = 0.
\]
The fluid is then excited by the action of external volume force supposed to capture ``in an equivalent manner'' the response
of the outer world, cf. Yakhot and Orszag \cite{YakOrs}. Such a hypothesis, when applied to realistic thermodynamically consistent models, however, gives rise to
a rather boring scenario: Either the system stabilizes to a static equilibrium, or the energy becomes infinite for the time
$t \to \infty$, see \cite{FP20}. This can be seen as a direct consequence of the Second law of thermodynamics as all mechanical energy is 
eventually converted to heat confined to the spatial domain by isolated boundary.  
To exhibit turbulent phenomena in the long run, the fluid system must be open; the coercive
effect of dissipation and thermal energy production being counterbalanced by the energy influx and outflux through the physical boundary. The aim of this work is to develop a mathematical theory for Newtonian models of compressible, viscous, and heat conducting fluid, with general in/out flow boundary conditions.

Motivated by Norman \cite{Norm}, we consider a bounded spatial domain $\Omega \subset R^d$ and suppose the fluid velocity
is given on $\partial \Omega$,
\begin{equation} \label{I1}
\vu|_{\partial \Omega} = \vuB.
\end{equation}
Furthermore, we decompose
\[
\partial \Omega = \Gamma_{\rm in} \cup \Gamma_{\rm out} \cup \Gamma_{\rm wall},
\]
\begin{equation} \label{i7}
\begin{split}
\Gamma_{{\rm in}} &= \left\{ x \in \partial \Omega \ \Big|\ \vc{n}(x) \cdot \vuB(x) < 0 \right\},\\
\Gamma_{{\rm out}} &= \left\{ x \in \partial \Omega \ \Big|\ \vc{n}(x) \cdot \vuB(x) > 0 \right\},\\
\Gamma_{{\rm wall}} &= \left\{ x \in \partial \Omega \ \Big|\ \vc{n}(x) \cdot \vuB(x) = 0 \right\},
\end{split}
\end{equation}
where $\vc{n}$ denotes the outer normal vector. The fluid mass density $\vr$ is given on the inflow boundary,
\begin{equation} \label{i8}
\vr = \vr_b \ \mbox{on}\ \Gamma_{\rm in}.
\end{equation}
Finally, denoting $e$ and $\vc{q}$ the internal energy and its diffusive (heat) flux, respectively, we prescribe
\begin{equation} \label{i9}
\Big[ \vr_b e \vuB  + \vc{q} \Big] \cdot \vc{n} = F_{i,b} \ \mbox{on}\ \Gamma_{\rm in},
\end{equation}
and
\begin{equation} \label{i10}
\vc{q} \cdot \vc{n} = 0 \ \mbox{on}\ \Gamma_{\rm wall} \cup \Gamma_{\rm out}.
\end{equation}
Here, $F_{i,b}$ is a given flux function reflecting the way the energy is flowing in/out of the physical domain.

{The boundary conditions (\ref{I1}--\ref{i10}) are adequate not only for the explanation of turbulent phenomena but they are also realistic in the modeling
of many real word applications. In fact, this is a natural and basic abstract setting for flows in pipelines, wind tunnels, turbines and jet engines, to name a few specific examples.}

This paper is the first attempt to establish the existence and weak-strong stability 
for the Navier-Stokes-Fourier system describing viscous compressible and heat conducting fluids. To the best of our knowledge, so far, all papers treating the
Navier-Stokes-Fourier equations in various settings deal always with periodic or homogenous boundary conditions for the velocity and
for the heat flux, see e.g. Bresch, Desjardins \cite{BrDe} and the results in \cite{FeNo6}, \cite{FeiNov10}, \cite{FeNo6A}, \cite{FeNoRoyal}, \cite{FENOSU}.

The mathematical theory developed in the present paper is based on the concept of weak (distributional) solutions, in the spirit
of Leray \cite{Ler} (incompressible fluids), Lions \cite{LI4} and \cite{FNP} (compressible barotropic fluids), and \cite{FeNo6A} (compressible and heat conducting fluids). The compressible Navier--Stokes system in the barotropic regime with inhomogeneous boundary conditions
\eqref{I1}, \eqref{i8} has been recently investigated in \cite{ChJiNo} (preceded by Girinon \cite{Girinon}) as far as existence of weak solutions is concerned,
and in \cite{AbFeNo}, \cite{KwNo} as far as the weak strong uniqueness is concerned. 
Similarly to \cite{FeNo6A}, our approach is based on careful implementation of the Second law of thermodynamics, in particular the existence of \emph{entropy} $s$ interrelated to the pressure $p$, the density $\vr$, the
internal energy $e$, and the (absolute) temperature $\vt$ through
\emph{Gibbs' equation}:
\begin{equation} \label{i2}
\vt D s = D e + p D \left( \frac{1}{\vr} \right).
\end{equation}
Besides a number of technical difficulties, the inhomogeneous boundary conditions require to control the state variables,
in particular the density, also on the outflow boundary, where their (normal) traces are interpreted in a very week sense.
Fortunately, the problem can be handled by convexity arguments on condition that the constitutive equations
satisfy the \emph{hypothesis of thermodynamic stability} specified below, cf. Bechtel, Rooney, Forest \cite{BeRoFo}.

\subsection{Field equations}

The motion of a general compressible viscous fluid is governed by the system of equations
\begin{equation} \label{i1}
\begin{split}
\partial_t \vr + \Div (\vr \vu) &= 0\\
\partial_t (\vr \vu) + \Div (\vr \vu \otimes \vu) + \Grad p &= \Div \mathbb{S} + \vr \vc{g}\\
\partial_t (\vr e) + \Div (\vr e \vu) + \Div \vc{q} &= \mathbb{S} : \Grad \vu - p \Div \vu,
\end{split}
\end{equation}
where $\vc{g}$ denotes the external driving force. We focus on linearly viscous fluids, where the viscous stress
tensor $\mathbb{S}$ is given by \emph{Newton's rheological law}
\begin{equation} \label{i3}
\mathbb{S}(\Grad \vu) = \mu \left( \Grad \vu + \Grad^t \vu - \frac{2}{d} \Div \vu \mathbb{I} \right) + \eta \Div \vu \mathbb{I},\
\mu > 0, \ \eta \geq 0.
\end{equation}
In addition, we impose \emph{Fourier's law}
\begin{equation} \label{i4}
\vc{q} = - \kappa \Grad \vt
\end{equation}
relating the heat flux to the temperature gradient. The system \eqref{i1}--\eqref{i4} is termed
\emph{Navier--Stokes--Fourier} system.

\subsection{First and Second law of thermodynamics}

The thermodynamic functions are interrelated through Gibbs' equation \eqref{i2}. In what follows, we alternatively consider
the \emph{standard} thermodynamic variables $(\vr, \vu, \vt)$ and the conservative--entropy variables $(\vr, \vm, S)$, where
$\vm = \vr \vu$ is the momentum, and $S = \vr s$ the total entropy. In particular, the pressure $p$ and internal energy $e$
may be viewed as $p = p(\vr, \vt)$ or $p=p(\vr, S)$, and, similarly, $e = e(\vr, \vt)$ or $e=e(\vr, S)$. To avoid confusion when partial derivatives are considered, we denote
\[
\frac{ \partial p(\vr, \vt) }{\partial \vr} = \frac{\partial p|_\vt}{\partial \vr},\
\frac{ \partial p(\vr, S) }{\partial \vr} = \frac{\partial p|_S}{\partial \vr},\
\frac{ \partial p(\vr, \vt) }{\partial \vt} = \frac{\partial p|_\vr}{\partial \vt},\
\frac{ \partial p(\vr, S) }{\partial S} = \frac{\partial p|_\vr}{\partial S}
\]
and similarly for $e$. The reason for using the standard variables is mainly because the
diffusive fluxes $\mathbb{S}$, $\vc{q}$ are easier to express in the standard variables, while
the conservative--entropy variables are more suitable in the weak formulation as they admit well defined traces, in particular the initial values, in the physical space time.

The \emph{thermodynamics stability hypothesis} written in terms of the standard variables reads
\begin{equation} \label{i5a}
\frac{\partial p|_\vt }{\partial \vr} > 0,\
\frac{\partial e|_\vr }{\partial \vt} > 0.
\end{equation}
The same condition may be expressed in the conservative--entropy variables as
\begin{equation} \label{i5}
E_{\rm int}(\vr, S) \equiv \vr e(\vr, S)\ \mbox{is a convex function of}\ (\vr, S), \ S = \vr s,
\end{equation}
see Section \ref{R}.
Moreover, it is straightforward to check that
\begin{equation} \label{i5b}
\frac{\partial (\vr e)|_S }{\partial \vr} = e - \vt s + \frac{p}{\vr},\
\frac{\partial (\vr e)|_\vr }{\partial S} = \vt,
\end{equation}
where the latter equality may be viewed as a definition of the temperature $\vt$ in the framework of the entropy--conservative variables.

It is easy to deduce from \eqref{i1} the \emph{energy equation}
\begin{equation} \label{i14}
\partial_t \left( \frac{1}{2} \vr {|\vu|^2} + \vr e \right)
+ \Div \left[ \left( \frac{1}{2} \vr {|\vu|^2} + \vr e + p \right) \vu \right] + \Div \vc{q}-
\Div \left( \mathbb{S} \cdot \vu \right) = \vr \vc{g} \cdot \vu,
\end{equation}
and, by virtue of Gibbs' relation \eqref{i2}, the \emph{entropy equation}
\begin{equation} \label{i11}
\partial_t (\vr s) + \Div (\vr s \vu) + \Div \left( \frac{\vc{q}}{\vt} \right) =
\frac{1}{\vt} \left( \mathbb{S} : \Grad \vu - \frac{\vc{q} \cdot \Grad \vt}{\vt} \right)
\end{equation}
Note that the equations \eqref{i14}, \eqref{i11}, and the internal energy equation $(\mbox{\ref{i1}})_3$ are equivalent in the framework
of regular solutions. In the weak formulation, the entropy balance is usually replaced by \emph{inequality}
\begin{equation} \label{i11a}
\partial_t (\vr s) + \Div (\vr s \vu) + \Div \left( \frac{\vc{q}}{\vt} \right) \geq
\frac{1}{\vt} \left( \mathbb{S} : \Grad \vu - \frac{\vc{q} \cdot \Grad \vt}{\vt} \right)
\end{equation}
The energy flux boundary condition \eqref{i9} can be expressed in terms of entropy as
\begin{equation} \label{i12}
\left[ \vr_b s(\vr_b, \vt) \vuB + \frac{\vc{q}}{\vt} \right] \cdot \vc{n} = S_{i,b} \ \mbox{on}\ \Gamma_{\rm in},
\end{equation}
where
\begin{equation} \label{i13}
S_{i,b} = \frac{F_{i,b}}{\vt} +  \Big[  s(\vr_b, \vt)  - \frac{e(\vr_b, \vt)}{\vt} \Big] \vr_b \vuB \cdot \vc{n}
\ \mbox{on}\ \Gamma_{\rm in}
\end{equation}

Finally, we recall the equation for the \emph{total energy}
\[
\intO{ E },\ E(\vr,\vu, e) = \frac{1}{2} \vr |\vu|^2 + \vr e.
\]
To this end, we first extend the boundary velocity $\vuB$ inside $\Omega$. After a straightforward manipulation, we deduce
\begin{equation} \label{i17}
\begin{split}
\frac{{\rm d}}{\dt} \intO{ \left( \frac{1}{2} \vr  |\vu - \vuB|^2  + \vr e \right) }
&+ \int_{\Gamma_{\rm out}} \vr e   \vuB \cdot \vc{n}
\ \D \sigma_x  \\ &= \intO{ \mathbb{S} : \Grad \vuB }   + \intO{ \vr \vc{g} \cdot (\vu - \vuB) } +
\frac{1}{2} \intO{ \vr \vu \cdot \Grad |\vuB|^2 }\\
&- \intO{ \Big( \vr \vu \otimes \vu + p \mathbb{I} \Big) : \Grad \vuB } -
\int_{\Gamma_{\rm in}} F_{i,b}
\ \D \sigma_x.
\end{split}
\end{equation}

\subsection{Mathematical theory in the framework of weak solutions}

The paper is organized as follows:

\begin{itemize}

\item In Section \ref{w}, we introduce the weak formulation of the problem. The leading idea is the same as in \cite{FeNo6A}, namely replacing the energy equation by the entropy inequality and the total energy balance. The completely new ingredient is suitable
accommodation of the boundary conditions. It turns out that the velocity $\vu$ as well as the temperature $\vt$ admit well defined
traces while the density $\vr$ does not. Moreover, it is convenient to include also the traces on the outflow part of the boundary to ensure stability of the solution set.

\item In Section \ref{R}, we derive a variant of the \emph{relative energy inequality} satisfied by \emph{any} weak solution
of the problem. The relative energy represents a Bregman distance (cf. e.g. Sprung \cite{Sprung}) between a weak solution and
an arbitrary trio of functions ranging in the associated phase space.

\item In Section \ref{W}, we show the weak--strong uniqueness principle. Any weak solution coincides with the strong solution
emanating from the sama initial/boundary data as long as the latter solution exists. The proof is an application of the relative energy
inequality.

\item Finally, in Section \ref{E}, we introduce an approximate scheme and prove existence of global--in--time weak solution for any physically admissible data.

\end{itemize}

\section{Weak formulation}
\label{w}

The weak formulation combines the ideas of \cite{FeNo6A} with those of \cite{ChJiNo} to accommodate the boundary data. We write down the field equations in terms of the standard variables $(\vr, \vu, \vt)$, however, the integrals on the outflow boundary
will be expressed in terms of the conservative--entropy variables $\vr$, $S$, and the internal energy $E_{\rm int}(\vr, S)$.
Accordingly, we shall always tacitly assume that any weak solution belongs at least to the class:
\begin{equation} \label{wc1}
\begin{split}
\vr &\in L^\infty(0,T; L^\gamma (\Omega)) \cap L^1((0,T) \times \Gamma_{\rm out}; \dt \times |\vuB \cdot \vc{n}| \dx)\ \mbox{for some}\ \gamma > 1, \\ &\vr \geq 0 \ \mbox{a.a. in}\ (0,T) \times \Omega;\\
\vu &\in L^q(0,T; W^{1,q}(\Omega; R^d)) \ \mbox{for some}\ q > 1,\
\vr \vu \in L^\infty(0,T; L^{\frac{2 \gamma}{\gamma + 1}} (\Omega; R^d));\\
\vt, \ \log(\vt) &\in L^2(0,T; W^{1,2}(\Omega)), \ \vt > 0 \ \mbox{a.a. in} \ (0,T) \times \Omega, \ 
\frac{1}{\vt} \in L^1((0,T) \times \Gamma_{\rm in}); \\ 
S &\in L^\infty(0,T; L^1(\Omega)) \cap L^1((0,T) \times \Gamma_{\rm out}),\
E_{\rm int}(\vr, S) \in L^1((0,T) \times \Gamma_{\rm out}).
\end{split}
\end{equation}

\begin{Definition}[Weak solution] \label{wD1}

Let $\Omega \subset R^d$, $d=2,3$ be a bounded domain with smooth boundary. Let the boundary data 
\[
\vuB \in C^2(\Ov{\Omega}; R^d), \ \vr_b \in C^1(\Ov{\Omega}; R^d), \ F_{i,b} \in C(\Ov{\Omega}),
\]
and the volume force
\[
\vc{g} \in C(\Ov{\Omega}; R^d)
\]
be given functions of $x \in \Ov{\Omega}$. 

We say that $(\vr, \vu, \vt)$ is a \emph{weak solution} to the Navier--Stokes--Fourier system
in $(0,T) \times \Omega$ if the following holds:

\begin{itemize}

\item
{\bf Equation of continuity}

\begin{equation} \label{w1}
\begin{split}
\left[ \intO{ \vr \varphi } \right]_{t = 0}^{t=\tau} &+ \int_0^\tau \int_{\Gamma_{\rm in}} \varphi \vr_b \vuB \cdot \vc{n} \
\D \sigma_x \dt + \int_0^\tau \int_{\Gamma_{\rm out}} \varphi \vr  \vuB \cdot \vc{n}\ 
\D \sigma_x \dt \\ &=
\int_0^\tau \intO{ \Big[ \vr \partial_t \varphi + \vr \vu \cdot \Grad \varphi \Big] } \dt
\end{split}
\end{equation}
holds for any $0 \leq \tau \leq T$,
and any $\varphi \in C^1([0,T] \times \Ov{\Omega})$;

\item
{\bf Momentum equation}

\begin{equation} \label{w2}
\begin{split}
\left[ \intO{ \vr \vu \cdot \bfphi } \right]_{t = 0}^{t=\tau}
&= \int_0^\tau \intO{ \Big[ \vr \vu \cdot \partial_t \bfphi + \vr \vu \otimes \vu : \Grad \bfphi + p(\vr,\vt) \Div \bfphi \Big] } \dt\\
&- \int_0^\tau \intO{ \mathbb{S} : \Grad \bfphi } \dt + \int_0^\tau \intO{ \vr \vc{g} \cdot \bfphi } \dt
\end{split}
\end{equation}
holds for any $0 \leq \tau \leq T$,
and any $\bfphi \in C^1_c([0,T] \times \Omega; R^d)$,
\begin{equation} \label{w3}
\vu - \vuB \in L^q(0,T; W^{1,q}_0(\Omega; R^d));
\end{equation}

\item
{\bf Total energy balance}

\begin{equation} \label{w4}
\begin{split}
&\left[ \intO{ \left( \frac{1}{2} \vr  |\vu - \vuB|^2  + \vr e \right) \psi } \right]_{t=0}^{t = \tau}
- \int_0^\tau \partial_t \psi\intO{ \left( \frac{1}{2} \vr  |\vu - \vuB|^2  + \vr e \right)  } \dt  \\
&+ \int_0^\tau \psi \int_{\Gamma_{\rm out}} E_{\rm int}(\vr, S)   \vuB \cdot \vc{n}
\ \D \sigma_x \dt \\ &\leq  \int_0^\tau \psi \intO{ \mathbb{S} : \Grad \vuB }\dt   + \int_0^\tau \psi \intO{ \vr \vc{g} \cdot (\vu - \vuB) }\dt +
\frac{1}{2} \int_0^\tau \psi \intO{ \vr \vu \cdot \Grad |\vuB|^2 } \dt \\
&- \int_0^\tau \psi \intO{ \Big( \vr \vu \otimes \vu + p \mathbb{I} \Big) : \Grad \vuB }\dt -
\int_0^\tau \psi \int_{\Gamma_{\rm in}} F_{i,b}
\ \D \sigma_x \dt
\end{split}
\end{equation}
holds for a.a. $0 \leq \tau \leq T$ and any $\psi \in C^1[0,T]$, $\psi \geq 0$;

\item
{\bf Entropy inequality}

\begin{equation} \label{w5}
\begin{split}
\left[ \intO{ \vr s \varphi } \right]_{t=0}^{t = \tau}
- \int_0^\tau &\intO{ \left[ \vr s \partial_t \varphi + \vr s \vu \cdot \Grad \varphi + \left( \frac{\vc{q}}{\vt} \right)
\cdot \Grad \varphi \right]} \dt
\\+ \int_0^\tau \int_{\Gamma_{\rm out}} \varphi S \vuB \cdot \vc{n} \ \D \sigma_x \dt
 &\geq
\int_0^\tau \intO{ \frac{\varphi}{\vt} \left( \mathbb{S} : \Grad \vu - \frac{\vc{q} \cdot \Grad \vt}{\vt} \right)} \dt \\
&- \int_0^\tau \int_{\Gamma_{\rm in}} \varphi \left( \frac{F_{i,b}}{\vt} +  \Big[  s(\vr_b, \vt)  - \frac{e(\vr_b, \vt)}{\vt} \Big] \vr_b \vuB \cdot \vc{n} \right)  \D \sigma_x \dt
\end{split}
\end{equation}
holds for a.a. $0 \leq \tau \leq T$, and
any $\varphi \in C^1([0,T] \times \Ov{\Omega})$, $\varphi \geq 0$.

\end{itemize}

\end{Definition}

The quantities $\vr \vuB \cdot \vc{n}|_{\Gamma_{\rm out}}$, $S \vuB \cdot \vc{n}|_{\Gamma_{\rm out}}$ can be (formally) identified with 
the normal traces of the fluxes $\vr \vu$, $\vr s \vu$, respectively, in the spirit of Chen, Torres, Ziemer \cite{ChToZi}. Their relation 
to the boundary integral containing $E_{\rm int}$ in \eqref{w4} is absolutely crucial for the property of stability of strong solutions in the class of weak solutions (weak--strong uniqueness principle).
The interested reader may consult \cite[Chapters 1--3]{FeNo6A} for a detailed discussion of the concept of weak solution introduced in
Definition \ref{wD1}. As we show in the next two sections, the weak solutions enjoy the important property of weak--strong uniqueness --
they coincide with the strong solution as long as the latter exists. To show this, however, certain technical hypotheses will be imposed on the constitutive relations.

\section{Relative energy as a Bregman distance}
\label{R}

The relative energy for the Navier--Stokes--Fourier system, written in the standard variables as
\[
E \left( \vr, \vu, \vt \ \Big| \tvr, \tvu, \tvt \right) = \frac{1}{2} \vr |\vu - \tvu|^2 +
H_{\tvt} (\vr, \vt) - \frac{ \partial H_{\tvt}(\tvr, \tvt)}{\partial \vr} (\vr - \tvr) - H_{\tvt}(\tvr, \tvt),
\]
\[
H_{\tvt}(\vr, \vt) \equiv \vr \Big( e(\vr, \vt) - \tvt s(\vr, \vt) \Big),
\]
was introduced in \cite{FeiNov10}. It is interesting to observe that
the relative energy represents a Bregman distance for the energy functional
\[
E(\vr, \vm, S) = \frac{1}{2} \frac{|\vm|^2}{\vr} + \vr e(\vr, S)
\]
written in terms of the conservative--entropy variables
\[
\vr, \ \vm, \ S = \vr s
\]
as long as the hypothesis of thermodynamics stability \eqref{i5} (or equivalently \eqref{i5a}) are satisfied.
Indeed it is easy to check, by virtue of \eqref{i5b},  
that
\[
E \left( \vr, \vm, S \ \Big| \tvr, \tvm, \tS \right)=E \left( \vr, \vu, \vt \ \Big| \tvr, \tvu, \tvt \right),\; \vm=\vr\vu,\,\tvm=\tvr\tvu,\, S=\vr e(\vr,\vt),\,\tS=\tvr e(\tvr,\tvt),
\]
where
\[
E \left( \vr, \vm, S \ \Big| \tvr, \tvm, \tS \right)= E(\vr, \vm, S) - \partial_{\vr, \vm, S} E(\tvr, \tvm, \tS) \cdot (\vr - \tvr, \vm - \tvm, S - \tS) -
E(\tvr, \tvm, \tS)
\]
as long as Gibbs' relation \eqref{i2} holds. As observed in \cite{FeiNov10}, the relative energy represent a distance between a potential weak solution $(\vr, \vu, \vt)$ and any trio of ``test functions'' $(\tvr, \tvu, \tvt)$. In particular, 
$E(\vr, \vm, S)$ is a convex function of the conservative--entropy variables, and the relative energy represents the associated 
Bregman distance.

More precisely, the mapping 
\begin{equation} \label{Edod}
\begin{split}
(\vr, \vt) &\mapsto (\vr, \vr s(\vr, \vt))\, \mbox{is a diffeomorphism}
\\
&\mbox{mapping}\ (0,\infty)^2 \ \mbox{onto an open convex set}\ \mathcal{E} \subset (0, \infty) \times R,
\end{split}
\end{equation} 
on which the internal energy 
\[
E_{\rm int} (\vr, S) = \vr e (\vr, S)
\]
is (strictly) convex. Extending 
\begin{equation} \label{Conv}
E_{\rm int} (\vr, S) = \left\{ \begin{array} {l} \infty \ \mbox{if} (\vr, S) \in \ R^2 \setminus \Ov{\mathcal{E}},\\ \\ 
\vr e(\vr, S) \ \mbox{if}\ (\vr, S) \in \mathcal{E},\\ \\
\liminf_{(\tvr, \widetilde{S}) \in \mathcal{E} ,(\tvr, \widetilde{S}) \to 
(\vr, S) } \tvr e(\tvr, \widetilde{S}) \ \mbox{if}\ (\vr, S) \in \partial \mathcal{E} \end{array}  \right.
\end{equation}
we obtain a convex l.s.c. function on $R^2$.


In the remaining part of this section, we derive a useful inequality satisfied by the relative energy if $(\vr, \vu, \vt)$ is a weak solution
of the Navier--Stokes--Fourier system.

\subsection{Derivation of the relative energy inequality}

We suppose that $(\tvr, \tvu, \tvt)$ are smooth functions of $(t,x) \in [0,T] \times \Ov{\Omega}$ satisfying the compatibility
condition
\begin{equation} \label{CC}
\tvu|_{\partial \Omega} = \vuB, \ 0 < \inf \tvr \leq \tvr \leq \sup \tvr < \infty,\ 0 < \inf \tvt \leq \tvt \leq \sup \tvt < \infty.
\end{equation}
{Starting from now we shall use abbreviated notation $\tilde e =e(\tvr,\tvt)$
and similarly for $\tilde p$, $\tilde s$, $\tilde S$, etc., whenever there is no danger of confusion.}

\subsubsection{Relative kinetic energy}

Consider $\tvu - \vuB$ as a test function in the momentum balance \eqref{w2}:
\[
\begin{split}
&\left[ \intO{ \vr \vu \cdot (\tvu - \vuB) } \right]_{t=0}^{t = \tau} \\ &=
\int_0^\tau \intO{ \Big[ \vr \vu \cdot \partial_t \tvu  + \vr \vu \otimes \vu : \Grad (\tvu - \vuB)
+ p \Div (\tvu - \vuB) - \mathbb{S} : \Grad (\tvu - \vuB) \Big] } \dt \\
&+ \int_0^\tau \intO{ \vr \vc{g} \cdot (\tvu - \vuB) } \dt
\end{split}
\]

Next, test the equation of continuity on $\frac{1}{2} (|\tvu|^2 - |\vuB|^2)$:
\[
\begin{split}
\left[ \intO{ \vr \left[ \frac{1}{2} \Big( |\tvu|^2 - |\vuB|^2 \Big) \right] } \right]_{t = 0}^{t = \tau}
=
\int_0^\tau \intO{ \Big[ \vr \partial_t \left(\frac{1}{2} |\tvu|^2 \right)  +
\vr \vu \cdot \Grad \frac{1}{2} \Big( |\tvu|^2 - |\vu_B|^2 \Big)  \Big] } \dt
\end{split}
\]

Finally, summing up the previous relations with the energy inequality \eqref{w4} yields
\[
\begin{split}
&\left[ \intO{ \left( \frac{1}{2} \vr  |\vu - \tvu|^2  + \vr e \right) } \right]_{t=0}^{t = \tau}
+ \int_0^\tau \int_{\Gamma_{\rm out}} E_{\rm int}(\vr, S)  \vuB \cdot \vc{n}
\ \D \sigma_x \dt 
\end{split}
\]
\begin{equation} \label{R1}
\begin{split}
 &\leq  \int_0^\tau \intO{ \mathbb{S} : \Grad \tvu }\dt   + \int_0^\tau \intO{ \vr \vc{g} \cdot (\vu - \tvu) }\dt +
\frac{1}{2} \int_0^\tau \intO{ \Big( \vr \partial_t |\tvu|^2 + \vr \vu \cdot \Grad |\tvu|^2 \Big) } \dt \\
&- \int_0^\tau \intO{ \vr \vu \cdot \partial_t \tvu } \dt - \int_0^\tau \intO{ \Big( \vr \vu \otimes \vu + p \mathbb{I} \Big) : \Grad \tvu }\dt -
\int_0^\tau \int_{\Gamma_{\rm in}} F_{i,b}
\ \D \sigma_x \dt
\end{split}
\end{equation}

\subsubsection{Entropy}

Recalling that
\[
\frac{ \partial (\vr e)|_\vr (\tvr, \tS) }{\partial S} = \tvt
\]
we use $\tvt$ as a test function in the entropy balance \eqref{w5} obtaining
\begin{equation} \label{R2}
\begin{split}
- \left[ \intO{ \vr s \tvt } \right]_{t=0}^{t=\tau}
+ \int_0^\tau &\intO{ \left[ \vr s \partial_t \tvt + \vr s \vu \cdot \Grad \tvt + \left( \frac{\vc{q}}{\vt} \right)
\cdot \Grad \tvt \right]} \dt
- \int_0^\tau \int_{\Gamma_{\rm out}} \tvt S \vuB \cdot \vc{n} \D \sigma_x \dt
\\ &\leq
- \int_0^\tau \intO{ \frac{\tvt}{\vt} \left( \mathbb{S} : \Grad \vu - \frac{\vc{q} \cdot \Grad \vt}{\vt} \right)} \dt \\
&+ \int_0^\tau \int_{\Gamma_{\rm in}} \left( F_{i,b} \frac{\tvt}{\vt} +  \tvt \Big[  s(\vr_b, \vt)  - \frac{e(\vr_b, \vt)}{\vt} \Big] \vr_b \vuB \cdot \vc{n} \right) \D \sigma_x \dt
\end{split}
\end{equation}

Summing up \eqref{R1}, \eqref{R2} and performing a simple manipulation  we obtain
\begin{equation} \label{R3}
\begin{split}
&\left[ \intO{ \left( \frac{1}{2} \vr  |\vu - \tvu|^2  + \vr e -
\frac{ \partial (\vr e)|_\vr (\tvr, \tS) }{\partial S} S - \tvr e(\tvr, \tS)
\right) } \right]_{t=0}^{t = \tau}\\&
+ \int_0^\tau \int_{\Gamma_{\rm out}} \left( E_{\rm int} (\vr, S) - \frac{ \partial (\vr e)|_\vr (\tvr, \tS) }{\partial S} S \right)  \vuB \cdot \vc{n}
\ \D \sigma_x \dt + \int_0^\tau \intO{ \frac{\tvt}{\vt} \left( \mathbb{S} : \Grad \vu - \frac{\vc{q} \cdot \Grad \vt}{\vt} \right)} \dt \\ &\leq \int_0^\tau \intO{ \Big( \vr (\tvu - \vu) \cdot \partial_t \tvu + \vr (\tvu - \vu) \otimes \vu : \Grad \tvu -
p \Div \tvu \Big)} \dt\\&+  \int_0^\tau \intO{ \mathbb{S} : \Grad \tvu }\dt   + \int_0^\tau \intO{ \vr \vc{g} \cdot (\vu - \tvu) }\dt  \\  &-  \int_0^\tau \intO{ \left[ \vr s \partial_t \tvt + \vr s \vu \cdot \Grad \tvt + \left( \frac{\vc{q}}{\vt} \right)
\cdot \Grad \tvt \right]} \dt   \\
&+ \int_0^\tau \int_{\Gamma_{\rm in}} \left( F_{i,b} \left( \frac{\tvt}{\vt} - 1 \right) +  \tvt \Big[  s(\vr_b, \vt)  - \frac{e(\vr_b, \vt)}{\vt} \Big] \vr_b \vuB \cdot \vc{n} \right) \D \sigma_x \dt
- \int_0^\tau \intO{ \partial_t (\tvr e(\tvr, \tS) ) } \dt.
\end{split}
\end{equation}

Next, testing the equation of continuity \eqref{w1} on $\frac{\partial (\vr e)|_S (\tvr, \tS) }{\partial \vr}$ we get
\[
\begin{split}
&\left[ \intO{ \vr \frac{\partial (\vr e)|_S (\tvr, \tS) }{\partial \vr} } \right]_{t = 0}^{t=\tau}  \\&+ \int_0^\tau \int_{\Gamma_{\rm in}} \vr_b \frac{\partial (\vr e)|_S (\tvr, \tS) }{\partial \vr}\vuB \cdot \vc{n}
\D \sigma_x \dt + \int_0^\tau \int_{\Gamma_{\rm out}} \vr \frac{\partial (\vr e)|_S (\tvr, \tS) }{\partial \vr} \vuB \cdot \vc{n}
\D \sigma_x \dt \\ &=
\int_0^\tau \intO{ \Big[ \vr \partial_t \frac{\partial (\vr e)|_S (\tvr, \tS) }{\partial \vr} + \vr \vu \cdot \Grad \frac{\partial (\vr e) (\tvr, \tS)|_S }{\partial \vr} \Big] } \dt
\end{split}
\]
Consequently,
\[
\begin{split}
&\left[ \intO{ \left( \frac{1}{2} \vr  |\vu - \tvu|^2  + \vr e -
\frac{ \partial (\vr e)|_\vr (\tvr, \tS) }{\partial S} S -
\frac{ \partial (\vr e)|_S (\tvr, \tS) }{\partial \vr} \vr - \tvr e(\tvr, \tS)
\right) } \right]_{t=0}^{t = \tau}\\&
+ \int_0^\tau \int_{\Gamma_{\rm out}} \left( E_{\rm int}(\vr, S) - \frac{ \partial (\vr e)|_\vr (\tvr, \tS) }{\partial S} S
- \frac{ \partial \vr e (\tvr, \tS)|_S }{\partial \vr} \vr \right)  \vuB \cdot \vc{n}
\ \D \sigma_x \dt \\&+ \int_0^\tau \intO{ \frac{\tvt}{\vt} \left( \mathbb{S} : \Grad \vu - \frac{\vc{q} \cdot \Grad \vt}{\vt} \right)} \dt 
\end{split}
\]
\[
\begin{split}
&\leq \int_0^\tau \intO{ \Big( \vr (\tvu - \vu) \cdot \partial_t \tvu + \vr (\tvu - \vu) \otimes \vu : \Grad \tvu -
p \Div \tvu \Big)} \dt\\&+  \int_0^\tau \intO{ \mathbb{S} : \Grad \tvu }\dt   + \int_0^\tau \intO{ \vr \vc{g} \cdot (\vu - \tvu) }\dt    \\&-  \int_0^\tau \intO{ \left[ \vr s \partial_t \tvt + \vr s \vu \cdot \Grad \tvt + \left( \frac{\vc{q}}{\vt} \right)
\cdot \Grad \tvt \right]} \dt   \\
&+ \int_0^\tau \int_{\Gamma_{\rm in}} \left( F_{i,b} \left( \frac{\tvt}{\vt} - 1 \right) +  \tvt \Big[  s(\vr_b, \vt)  - \frac{e(\vr_b, \vt)}{\vt} \Big] \vr_b \vuB \cdot \vc{n} \right) \D \sigma_x \dt\\ &+ \int_0^\tau \int_{\Gamma_{\rm in}} \vr_b \frac{\partial (\vr e)|_S (\tvr, \tS) }{\partial \vr}\vuB \cdot \vc{n}
\D \sigma_x \dt
- \int_0^\tau \intO{ \partial_t (\tvr e(\tvr, \tS) ) } \dt\\
&- \int_0^\tau \intO{ \Big[ \vr \partial_t \frac{\partial (\vr e)|_S (\tvr, \tS) }{\partial \vr} + \vr \vu \cdot \Grad \frac{\partial (\vr e)|_S (\tvr, \tS) }{\partial \vr} \Big] } \dt
\end{split}
\]

\subsubsection{Final form}

After a simple manipulation based on Gibbs' relation we deduce the final form of the relative energy inequality:\footnote{In what follows, we denote $\tilde e=e(\tilde\vr,\tilde\vt)$, $\tilde p=p(\tilde\vr,\tilde\vt)$ etc. if there is no
danger of confusion.}

\begin{equation} \label{R4}
\begin{split}
&\left[ \intO{ \left( \frac{1}{2} \vr  |\vu - \tvu|^2  + \vr e -
\frac{ \partial (\vr e)|_\vr (\tvr, \tS) }{\partial S} (S - \tS) -
\frac{ \partial (\vr e)|_S (\tvr, \tS) }{\partial \vr} (\vr - \tvr) - \tvr e(\tvr, \tS)
\right) } \right]_{t=0}^{t = \tau}\\&
+ \int_0^\tau \int_{\Gamma_{\rm out}} \left( E_{\rm int}(\vr, S) - \frac{ \partial (\vr e)|_\vr  (\tvr, \tS) }{\partial S} S
- \frac{ \partial (\vr e)_S (\tvr, \tS) }{\partial \vr} \vr \right)  \vuB \cdot \vc{n}
\ \D \sigma_x \dt \\&+ \int_0^\tau \intO{ \frac{\tvt}{\vt} \left( \mathbb{S} : \Grad \vu - \frac{\vc{q} \cdot \Grad \vt}{\vt} \right)} \dt 
\end{split}
\end{equation}
\[
\begin{split}
&\leq - \int_0^\tau \intO{ \vr (\tvu - \vu) \otimes (\tvu - \vu) : \Grad \tvu } \dt
 \\
&+ \int_0^\tau \intO{ \vr (\tvu - \vu) \cdot \Big( \partial_t \tvu + \tvu \cdot \Grad \tvu + \frac{1}{\tvr} \Grad \widetilde{p} - \vc{g} \Big)} \dt\\&- \int_0^\tau \intO{ p \Div \tvu } \dt + \int_0^\tau \intO{ \frac{\vr}{\tvr} (\vu - \tvu) \cdot \Grad \widetilde{p} } \dt +  \int_0^\tau \intO{ \mathbb{S} : \Grad \tvu }\dt       \\&-  \int_0^\tau \intO{ \left[ \vr ( s - \widetilde{s} ) \partial_t \tvt + \vr (s - \widetilde{s}) \vu \cdot \Grad \tvt + \left( \frac{\vc{q}}{\vt} \right)
\cdot \Grad \tvt \right]} \dt   \\
&+ \int_0^\tau \int_{\Gamma_{\rm in}} \left( F_{i,b} \left( \frac{\tvt}{\vt} - 1 \right) +  \tvt \Big[  s(\vr_b, \vt) -
\widetilde{s}  + \frac{\widetilde{e}}{\tvt} - \frac{e(\vr_b, \vt)}{\vt} \Big] \vr_b \vuB \cdot \vc{n} \right) \D \sigma_x \dt\\ &+ \int_0^\tau \int_{\Gamma_{\rm in}} \frac{\vr_b}{\tvr} \widetilde{p} \vuB \cdot \vc{n}
\D \sigma_x \dt
+ \int_0^\tau \intO{\left( \left( 1 - \frac{\vr}{\tvr} \right) \partial_t \widetilde{p} - \frac{\vr}{\tvr} \vu \cdot \Grad \widetilde{p} \right) } \dt .
\end{split}
\]
Here, 
for the sake of brevity, we have used the notation $\widetilde{b} = b(\tvr, \tvt)$. It is worth noting that the relative energy inequality \eqref{R4} coincides, modulo the boundary terms, with that obtained in \cite[Section 3.2] {EF101}.

We have shown the following result.

\begin{Proposition} [Relative energy inequality] \label{wP1}
Let a trio of continuously differentiable functions $(\tvr, \tvu, \tvt)$ belong to the class \eqref{CC}.

Then the relative energy inequality \eqref{R4} holds for any weak solution $(\vr, \vu, \vt)$ of the Navier--Stokes--Fourier
system in the sense of Definition \ref{wD1}.

\end{Proposition}

\section{Weak--strong uniqueness}
\label{W}

We now suppose that $(\tvr, \tvu, \tvt)$ is a regular solution {of the system 
(\ref{i1}--\ref{i4}) satisfying the boundary conditions (\ref{I1}--\ref{i10})} and belonging to the class \eqref{CC} and use it as a test function in \eqref{R4}.

\subsection{Momentum balance}

As $(\tvr, \tvu, \tvt)$ satisfies the momentum balance, we get
\[
\begin{split}
&\intO{ \vr (\tvu - \vu) \cdot \Big( \partial_t \tvu + \tvu \cdot \Grad \tvu + \frac{1}{\tvr} \Grad \widetilde{p} - \vc{g} \Big)}
+ \intO{ \mathbb{S} : \Grad \tvu } \\&=
\intO{ \frac{\vr}{\tvr} (\tvu - \vu) \cdot \Div \widetilde{\mathbb{S}} + \mathbb{S} : \Grad \tvu }\\
&=  \intO{ \left( \frac{\vr}{\tvr} - 1\right) (\tvu - \vu) \cdot \Div \widetilde{\mathbb{S}}}
+ \intO{ \Grad (\vu - \tvu) : \widetilde{\mathbb{S}} + \mathbb{S} : \Grad \tvu }.
\end{split}
\]
Consequently, after a straightforward manipulation, the relative energy inequality
\eqref{R4} gives rise to
\[
\begin{split}
&\left[ \intO{ \left( \frac{1}{2} \vr  |\vu - \tvu|^2  + \vr e -
\frac{ \partial (\vr e)|_\vr (\tvr, \tS) }{\partial S} (S - \tS) -
\frac{ \partial (\vr e)|_S (\tvr, \tS) }{\partial \vr} (\vr - \tvr) - \tvr e(\tvr, \tS)
\right) } \right]_{t=0}^{t = \tau}\\&
+ \int_0^\tau \int_{\Gamma_{\rm out}} \left( E_{\rm int}(\vr, S) - \frac{ \partial (\vr e)|_\vr (\tvr, \tS) }{\partial S} S
- \frac{ \partial (\vr e)|_S (\tvr, \tS) }{\partial \vr} \vr \right)  \vuB \cdot \vc{n}
\ \D \sigma_x \dt \\&+ \int_0^\tau \intO{ \left( \left( \frac{\tvt}{\vt} - 1 \right) \mathbb{S} : \Grad \vu -
\left(1 - \frac{\tvt}{\vt} \right) \frac{\vc{q} \cdot \Grad \vt}{\vt} \right)} \dt + \int_0^\tau \intO{
\frac{\vc{q}}{\vt} \cdot \left( \Grad \tvt - \Grad \vt \right) } \dt
\\
&+ \int_0^\tau \intO{ (\mathbb{S} - \widetilde{\mathbb{S}}): (\Grad \vu - \Grad \tvu) } \dt
\end{split}
\]
\begin{equation} \label{R5}
\begin{split}
&\leq - \int_0^\tau \intO{ p \Div \tvu } \dt + \int_0^\tau \intO{  (\vu - \tvu) \cdot \Grad \widetilde{p} } \dt   \\&
+\int_0^\tau \intO{ \left[ \tvr ( \widetilde{s} - s ) \left( \partial_t \tvt + \tvu \cdot \Grad \tvt \right)  \right]} \dt   \\
&+ \int_0^\tau \int_{\Gamma_{\rm in}} \left( F_{i,b} \left( \frac{\tvt}{\vt} - 1 \right) +  \tvt \Big[  s(\vr_b, \vt) -
\widetilde{s}  + \frac{\widetilde{e}}{\tvt} - \frac{e(\vr_b, \vt)}{\vt} \Big] \vr_b \vuB \cdot \vc{n} \right) \D \sigma_x \dt\\ &+ \int_0^\tau \int_{\Gamma_{\rm in}} \frac{\vr_b}{\tvr} \widetilde{p} \vuB \cdot \vc{n} \ 
\D \sigma_x \dt 
+ \int_0^\tau \intO{\left( \left( 1 - \frac{\vr}{\tvr} \right) \partial_t \widetilde{p} - \frac{\vr}{\tvr} \vu \cdot \Grad \widetilde{p} \right) } \dt 
+ \int_0^\tau {\rm Er}_1 (t) \dt,
\end{split}
\end{equation}
with an ``error term'' 
\begin{equation} \label{R5a}
\begin{split}
{\rm Er}_1 &= - \intO{ \vr (\tvu - \vu) \otimes (\tvu - \vu) : \Grad \tvu } 
+ \intO{\left(\frac{\tvr - \vr}{\tvr} \right) (\vu - \tvu) \cdot \left( \Grad \widetilde{p} - \Div \widetilde{
\mathbb{S} } \right) } 
\\
&+ \intO{ \vr (\widetilde{s} - s) (\vu - \tvu) \cdot \Grad \tvt }  +
\intO{ (\vr - \tvr) ( \widetilde{s} - s ) \left( \partial_t \tvt + \tvu \cdot \Grad \tvt \right)}. 
\end{split}
\end{equation}

\subsection{Pressure}

First observe that
\[
- \intO{ \tvu \cdot \Grad \widetilde{p} } = - \int_{\partial \Omega} \widetilde{p} \vuB \cdot \vc{n} \D \sigma_x +
\intO{ \widetilde{p} \Div \tvu }.
\]
Next, a direct manipulation yields 
\[
\widetilde{p} = \frac{ \partial (\vr e)|_S (\tvr, \tS) } \tS + \frac{ \partial (\vr e)|_\vr (\tvr, \tS) }{\partial \vr} \tvr - \tvr \widetilde{e}. 
\]
Finally, we report the identity 
\[
\begin{split}
&\left(1 - \frac{\vr}{\tvr} \right) \Big(\partial_t \widetilde{p} +  \tvu \cdot \Grad \widetilde{p} \Big)
+ \Big( \widetilde{p} - p \Big) \Div \tvu 
= \Div \tvu \left( \widetilde{p} - \frac{\partial \widetilde{p} }{\partial \vr} (\tvr - \vr) - \frac{\partial \widetilde{p} }{
\partial \vt} (\tvt - \vt) - p \right)\\
&- \tvr (\partial_t \tvt + \tvu \cdot \Grad \tvt ) \left( \frac{\partial \widetilde{s} }{\partial \vr}(\tvr - \vr) +
\frac{\partial \widetilde{s} }{\partial \vt}(\tvt - \vt) \right)\\
&+ \left(1 - \frac{\vt}{\tvt} \right) \left( \widetilde{\mathbb{S}}: \Grad \tvu - \frac{\widetilde{\vc{q}} \cdot \Grad \tvt}{\tvt} \right) + (\vt - \tvt) \Div \left( \frac{ \widetilde{\vc{q}} }{\tvt } \right),
\end{split}
\]
see \cite[Section 6]{EF101}.

Consequently, plugging these three relations in \eqref{R5} and using the fact $\tvr|_{\Gamma_{\rm in}} = \vr_b$, we may infer that
\begin{equation}\label{R6}
\begin{split}
&\left[ \intO{ \left( \frac{1}{2} \vr  |\vu - \tvu|^2  + \vr e -
\frac{ \partial (\vr e)|_\vr (\tvr, \tS) }{\partial S} (S - \tS) -
\frac{ \partial (\vr e)|_S (\tvr, \tS) }{\partial \vr} (\vr - \tvr) - \tvr e(\tvr, \tS)
\right) } \right]_{t=0}^{t = \tau}\\&
+ \int_0^\tau \int_{\Gamma_{\rm out}} \left( E_{\rm int}(\vr, S) - \frac{ \partial (\vr e)|_\vr (\tvr, \tS) }{\partial S} (S - \tS)
- \frac{ \partial (\vr e)|_S (\tvr, \tS) }{\partial \vr} (\vr - \tvr)  -  E_{\rm int}(\tvr, \tS) \right)  \vuB \cdot \vc{n}
\ \D \sigma_x \dt \\&+ \int_0^\tau \intO{ \left( \frac{\tvt}{\vt} - 1 \right) \mathbb{S} : \Grad \vu } \dt + \int_0^\tau \intO{ \left(\frac{\vt}{\tvt} - 1 \right) \left( \widetilde{\mathbb{S}}: \Grad \tvu \right) } \dt\\
& + \int_0^\tau \intO{ \left(1 - \frac{\vt}{\tvt} \right)\frac{\widetilde{\vc{q}} \cdot \Grad \tvt}{\tvt} } \dt + 
\int_0^\tau \intO{ \left(\frac{\tvt}{\vt} - 1 \right) \frac{\vc{q} \cdot \Grad \vt}{\vt} } \dt
\\
&+ \int_0^\tau \intO{ (\mathbb{S} - \widetilde{\mathbb{S}}): (\Grad \vu - \Grad \tvu) } \dt + \int_0^\tau \intO{
\left( \frac{\vc{q}}{\vt} - \frac{\widetilde{\vc{q}}}{\tvt} \right) \cdot \left( \Grad \tvt - \Grad \vt \right) } \dt
\\&\leq 
\int_0^\tau \int_{\Gamma_{\rm in}} \left( F_{i,b} \left( \frac{\tvt}{\vt} - 1 \right) +  \tvt \Big[  s(\vr_b, \vt) -
\widetilde{s}  + \frac{\widetilde{e}}{\tvt} - \frac{e(\vr_b, \vt)}{\vt} \Big] \vr_b \vuB \cdot \vc{n}
+ \widetilde{\vc{q}} \cdot \vc{n} \left( \frac{\vt}{\tvt} - 1 \right)
\right)
\D \sigma_x \dt \\ &+ \int_0^\tau \intO{ {\rm Er}_2(t) } \dt,
\end{split}
\end{equation}
with 
\begin{equation} \label{R6a}
\begin{split}
{\rm Er}_2 = &- \intO{ \vr (\tvu - \vu) \otimes (\tvu - \vu) : \Grad \tvu } 
+ \intO{\left(1 - \frac{\vr}{\tvr} \right) (\vu - \tvu) \cdot \left( \Grad \widetilde{p} - \Div \widetilde{
\mathbb{S} } \right) } 
\\
&+ \intO{ \vr (\widetilde{s} - s) (\vu - \tvu) \cdot \Grad \tvt }  +
\intO{ (\vr - \tvr) ( \widetilde{s} - s ) \left( \partial_t \tvt + \tvu \cdot \Grad \tvt \right)} 
 \\
&+
\intO{ \left(1 - \frac{\vr}{\tvr} \right) (\vu - \tvu) \cdot \Grad \widetilde{p} } \\&+ 
\intO{ \tvr \left( \widetilde{s} - \frac{\partial \widetilde{s}|_\vt }{\partial \vr}(\tvr - \vr) -
\frac{\partial \widetilde{s}|_\vr  }{\partial \vt} (\tvt - \vt)- s \right) \left( \partial_t \tvt + \tvu \cdot \Grad \tvt \right)  }
\end{split}
\end{equation}

\subsection{Conclusion}

As a consequence of the hypothesis of thermodynamic stability expressed via convexity of the function $E_{\rm int}(\vr,S)$, we get 
\[
\int_0^\tau \int_{\Gamma_{\rm out}} \left( E_{\rm int}(\vr, S) - \frac{ \partial (\vr e)|_\vr (\tvr, \tS) }{\partial S} (S - \tS)
- \frac{ \partial (\vr e)|_S (\tvr, \tS) }{\partial \vr} (\vr - \tvr)  -  E_{\rm int}(\tvr, \tS) \right)  \vuB \cdot \vc{n}
\ \D \sigma_x \dt \geq 0.
\]
Consequently, inequality \eqref{R6} can be rewritten in terms of the standard variables as 
\begin{equation}\label{R7}
\begin{split}
&\left[ \intO{ E \left( \vr, \vu, \vt \ \Big| \tvr, \tvu, \tvt \right) } \right]_{t=0}^{t = \tau}\\&+ \int_0^\tau \intO{ \left( \frac{\tvt}{\vt} - 1 \right) \mathbb{S} : \Grad \vu } \dt + \int_0^\tau \intO{ \left(\frac{\vt}{\tvt} - 1 \right) \left( \widetilde{\mathbb{S}}: \Grad \tvu \right) } \dt\\
& + \int_0^\tau \intO{ \left(1 - \frac{\vt}{\tvt} \right)\frac{\widetilde{\vc{q}} \cdot \Grad \tvt}{\tvt} } \dt + 
\int_0^\tau \intO{ \left(\frac{\tvt}{\vt} - 1 \right) \frac{\vc{q} \cdot \Grad \vt}{\vt} } \dt
\\
&+ \int_0^\tau \intO{ (\mathbb{S} - \widetilde{\mathbb{S}}): (\Grad \vu - \Grad \tvu) } \dt + \int_0^\tau \intO{
\left( \frac{\vc{q}}{\vt} - \frac{\widetilde{\vc{q}}}{\tvt} \right) \cdot \left( \Grad \tvt - \Grad \vt \right) } \dt
\\&\leq 
\int_0^\tau \int_{\Gamma_{\rm in}} \left( F_{i,b} \left( \frac{\tvt}{\vt} - 1 \right) +  \tvt \Big[  s(\vr_b, \vt) -
\widetilde{s}  + \frac{\widetilde{e}}{\tvt} - \frac{e(\vr_b, \vt)}{\vt} \Big] \vr_b \vuB \cdot \vc{n}
+ \widetilde{\vc{q}} \cdot \vc{n} \left( \frac{\vt}{\tvt} - 1 \right)
\right)
\D \sigma_x \dt \\ &+ \int_0^\tau \intO{ {\rm Er}_2(t) } \dt.
\end{split}
\end{equation}
With the exception of the boundary integral, this is the same inequality as in \cite[Section 6, formula (71)]{EF101}. 
We may therefore anticipate results similar to \cite{EF101} as soon as we handle the boundary terms.  

\subsubsection{Inflow boundary}

Since $(\tvr, \tvu, \tvt)$ is a strong solution, 
\[
\vr_b \widetilde{e} \vuB \cdot \vc{n} + \widetilde{\vc{q}} \cdot \vc{n} = F_{i,b},
\]
and the boundary integral reads
\[
\int_{\Gamma_{\rm in}} \left( F_{i,b} \left( \frac{\vt}{\tvt} + \frac{\tvt}{\vt} - 2 \right) +  \tvt \Big[  s(\vr_b, \vt) -
\widetilde{s}  + \frac{\widetilde{e}}{\tvt} - \frac{e(\vr_b, \vt)}{\vt} \Big] \vr_b \vuB \cdot \vc{n}
- \widetilde{e} \left( \frac{\vt}{\tvt} - 1 \right) \vr_b \vuB \cdot \vc{n}
\right)
\D \sigma_x
\]

Next, it follows from Gibbs' equation that
\[
\widetilde{e} - \vt \widetilde{s} =
e(\vr_b, \tvt) - \vt s(\vr_b, \tvt) \geq e(\vr_b, \vt) - \vt s(\vr_b, \vt).
\]
Indeed, as shown in \cite{FeiNov10}, the function 
\[
\tvt \mapsto e(\vr_b, \tvt) - \vt s(\vr_b, \tvt) 
\]
is non--negative attaining its minimum at $\tvt = \vt$ for any fixed $\vr_b$, $\vt$. 
Consequently,
\[
\begin{split}
\tvt \Big[  &s(\vr_b, \vt) -
\widetilde{s}  + \frac{\widetilde{e}}{\tvt} - \frac{e(\vr_b, \vt)}{\vt} \Big]
= \frac{\tvt}{\vt} \Big[ \vt s(\vr_b, \vt) -
\vt \widetilde{s}  + \frac{\vt \widetilde{e}}{\tvt} - {e(\vr_b, \vt)} \Big] \\ &=
\frac{\tvt}{\vt} \Big[ \vt s(\vr_b, \vt) -
\vt \widetilde{s}  + \widetilde{e} - {e(\vr_b, \vt)} \Big] + \widetilde{e} - \frac{\tvt}{\vt} \widetilde{e};
\end{split}
\]
whence the boundary integral can be controlled as 
\begin{equation} \label{R7a}
\begin{split}
\int_{\Gamma_{\rm in}} &\left( F_{i,b} \left( \frac{\vt}{\tvt} + \frac{\tvt}{\vt} - 2 \right)
- \widetilde{e} \left( \frac{\vt}{\tvt} + \frac{\tvt}{\vt} - 2 \right) \vr_b \vuB \cdot \vc{n}
\right)
\D \sigma_x \\
&+ \int_{\Gamma_{\rm in}} \frac{\tvt}{\vt} \Big[ \vt s(\vr_b, \vt) -
\vt \widetilde{s}  + \widetilde{e} - {e(\vr_b, \vt)} \Big] \vr_b \vuB \cdot \vc{n} \D \sigma_x \leq \int_{\Gamma_{\rm in}} \left( \frac{\vt}{\tvt} + \frac{\tvt}{\vt} - 2 \right) \widetilde{\vc{q}} \cdot \vc{n}\ 
\D \sigma_x.
\end{split}
\end{equation}

\subsubsection{Conditional weak--strong uniqueness}

We suppose that the weak solution $(\vr, \vu, \vt)$ belongs to the ``non--degenerate'' area 
\begin{equation} \label{R8}
0 < \underline{\vr} \leq \vr(t,x) \leq \Ov{\vr}, \ 0 < \underline{\vt} \leq \vt(t,x) \leq \Ov{\vt} \ \mbox{for a.a.}\ (t,x) \in 
(0,T) \times \Omega, 
\end{equation}
where $\underline{\vr}, \underline{\vt}, \Ov{\vr}, \Ov{\vt}$ are constants. Under these circumstances, we can 
apply a Gronwall type argument to the inequality \eqref{R7} exactly  
as in \cite[Section 6.1, Theorem 6.1]{EF101} to show a conditional weak--strong uniqueness result. The boundary integral 
\eqref{R7a} can be handled by means of the trace theorem and interpolation as 
\begin{equation} \label{R9} 
\begin{split}
\int_{\Gamma_{\rm in}} \left( \frac{\vt}{\tvt} + \frac{\tvt}{\vt} - 2 \right) \widetilde{\vc{q}} \cdot \vc{n}
\D \sigma_x &\aleq \int_{\Gamma_{\rm in}} |\vt - \tvt|^2 \ \D \sigma_x \aleq \left\| \vt - \tvt \right\|_{W^{\alpha,2}(\Omega)}^2\\
&\leq \delta \| \vt - \tvt \|^2_{W^{1,2}(\Omega)} + c(\delta) \| \vt - \tvt \|^2_{L^2(\Omega)},\ \frac{1}{2} < \alpha < 1,
\end{split}
\end{equation}
where $\delta > 0$ can be chosen arbitrarily small. Thus the same arguments as in \cite[Section 6.1]{EF101} 
give rise to the following result. 

\begin{Theorem}[Conditional weak--strong uniqueness] \label{RT1}

Let the thermodynamic functions $p$, $e$, and $s$ be continuously differentiable functions of $\vr$ and $\vt$ satisfying 
Gibbs' equation \eqref{i1}, together with the hypothesis of thermodynamics stability {\eqref{i5}} 
(specified in \eqref{Edod}, \eqref{Conv}).
Let the transport 
coefficients $\mu$, $\eta$, $\kappa$ be continuously differentiable functions of $\vr$ and $\vt$, 
\[
\mu > 0,\ \kappa > 0, \eta \geq 0.
\]
Let $(\vr, \vu, \vt)$ be a weak solution
of the Navier--Stokes system \eqref{i1}--\eqref{i4}, with the boundary conditions \eqref{I1}--\eqref{i10},  
in the sense specified in Definition \ref{wD1} satisfying 
\[
0 < \underline{\vr} \leq \vr(t,x) \leq \Ov{\vr}, \ 0 < \underline{\vt} \leq \vt(t,x) \leq \Ov{\vt} \ \mbox{for a.a.}\ (t,x) \in 
(0,T) \times \Omega. 
\]
Suppose that the same problem (with the same initial and boundary data) admits a strong solution $(\tvr, \tvu, \tvt)$ 
in the class 
\[
\tvr, \ \tvu, \ \tvt \in 
C^1([0,T] \times \Ov{\Omega}), \ \partial^2_x \tvu, \ \partial^2_x  \tvt \in 
C([0,T] \times \Ov{\Omega}). 
\]

Then 
\[
\vr = \tvr,\ \vu = \tvu, \ \vt = \tvt \ \mbox{in}\ [0,T] \times \Ov{\Omega}.
\]

\end{Theorem}

\subsubsection{Unconditional weak--strong uniqueness}
\label{UWU}

Unfortunately, the existence result proved below does not provide weak solutions ranging in the physically 
``regular'' domain \eqref{R8}. To save the weak--strong uniqueness principle, we are forced to impose certain technical hypotheses 
on the constitutive relations. Motivated by \cite[Chapters 2,3]{FeNo6A}, we suppose that 
the pressure $p$ obeys a state equation in the form
\begin{equation} \label{ws1}
p(\vr, \vt) = \vt^{5/2} P \left( \frac{ \vr }{\vt^{3/2}} \right) + \frac{a}{3} \vt^4,\ a > 0,
\end{equation}
with $P \in C^1[0,\infty)$. 
In accordance with Gibbs' equation (\ref{i2}), we get
\begin{equation} \label{ws2}
e(\vr, \vt) = \frac{3}{2} \frac{\vt^{5/2}}{\vr} P \left( \frac{ \vr }{\vt^{3/2}} \right) + \frac{a}{\vr} \vt^4,
\end{equation}
and
\begin{equation} \label{ws3}
s(\vr, \vt) = \mathcal{S} \left( \frac{\vr} {\vt^{3/2}} \right) + \frac{4a}{3} \frac{\vt^3}{\vr},
\end{equation}
where
\begin{equation} \label{ws4}
\mathcal{S}'(Z) = - \frac{3}{2} \frac{ \frac{5}{3} P(Z) - P'(Z) Z }{Z^2}.
\end{equation}
Moreover, the thermodynamic stability requires
\begin{equation} \label{ws5}
P'(Z) > 0 \ \mbox{for any}\ Z \geq 0,\ \frac{ \frac{5}{3} P(Z) - P'(Z) Z }{Z} > 0 \ \mbox{for any}\ Z > 0.
\end{equation}
In particular,  the function $Z \mapsto P(Z) / Z^{5/3}$ is decreasing, and we suppose
\begin{equation} \label{ws6}
\lim_{Z \to \infty} \frac{ P(Z) }{Z^{5/3}} = p_\infty > 0.
\end{equation}
Next, we impose technical but physically grounded hypothesis (see \cite[Chapter 2]{FeNo6A})
\begin{equation} \label{ws7}
P(0) = 0,\ \frac{ \frac{5}{3} P(Z) - P'(Z) Z }{Z} < c \ \mbox{for all}\ Z > 0.
\end{equation}

Finally, we require the transport coefficients to be continuously differentiable for $\vt \in [0,\infty)$,
\begin{equation} \label{ws8}
\underline{\mu} (1 + \vt^\Lambda) \leq \mu(\vt) \leq \Ov{\mu} (1 + \vt^\Lambda),\ | \mu'(\vt) | < c \ \mbox{for all}\ \vt \in [0, \infty),
\ \frac{2}{5} < \Lambda \leq 1,
\end{equation}
\begin{equation} \label{ws9}
0 \leq \eta (\vt) \leq \Ov{\eta} (1 + \vt^\Lambda) \ \mbox{for all}\ \vt \in [0, \infty),
\end{equation}
\begin{equation} \label{ws10}
\underline{\kappa} (1 + \vt^3) \leq \kappa (\vt) \leq \Ov{\kappa}(1 + \vt^3) \ \mbox{for all}\ \vt \in [0, \infty).
\end{equation}

{Observe that function ${\cal S}$ is decreasing on $(0,\infty)$ and we can suppose without loss of generality $\lim_{Z\to\infty}{\cal S}(Z)\in \{-\infty,0\}$. If 
\begin{equation}\label{dod4}
\lim_{Z\to\infty} {\cal S}(Z)=0
\end{equation}
then $s$ satisfies the Third law of thermodynamics, {cf.
Belgiorno \cite{BEL1}, \cite{BEL2}.}
The reader may consult \cite[Chapter 1]{FeNo6A} for the physical background of the above hypotheses.

The important observation made in \cite[Section 6.1, formula (78)]{EF101} is the following inequality: 
\begin{equation} \label{ws11}
\begin{split}
\left\| \vt - \tvt \right\|^2_{W^{1,2}(\Omega)}  \aleq
&\intO{ \left(1 - \frac{\vt}{\tvt} \right)\frac{\widetilde{\vc{q}} \cdot \Grad \tvt}{\tvt} } + 
\intO{ \left(\frac{\tvt}{\vt} - 1 \right) \frac{\vc{q} \cdot \Grad \vt}{\vt} }
\\
& + \intO{
\left( \frac{\vc{q}}{\vt} - \frac{\widetilde{\vc{q}}}{\tvt} \right) \cdot \left( \Grad \tvt - \Grad \vt \right) } 
+ \intO{ E \left( \vr, \vu, \vt \ \Big| \ \tvr, \tvu, \tvt \right) } 
\end{split}
\end{equation}

The bound \eqref{ws11} allows us to control the boundary integral 
\[
\int_{\Gamma_{\rm in}} \left( F_{i,b} \left( \frac{\vt}{\tvt} + \frac{\tvt}{\vt} - 2 \right) +  \tvt \Big[  s(\vr_b, \vt) -
\widetilde{s}  + \frac{\widetilde{e}}{\tvt} - \frac{e(\vr_b, \vt)}{\vt} \Big] \vr_b \vuB \cdot \vc{n}
- \widetilde{e} \left( \frac{\vt}{\tvt} - 1 \right) \vr_b \vuB \cdot \vc{n}
\right)
\D \sigma_x
\]
as long as further restrictions are imposed relating $F_{i,b}$, $\vr_b$, and the structural constant $p_\infty$ in 
\eqref{ws6}. Following Norman \cite[formula (2.10b)]{Norm} we suppose that
\begin{equation} \label{ws12}
F_{i,b}(x) < 0 \ \mbox{on}\ \Gamma_{\rm in}.
\end{equation}
Our goal is to show, similarly to \eqref{R9}, that  
\begin{equation} \label{ws13}
\begin{split}
\int_{\Gamma_{\rm in}} &\left( F_{i,b} \left( \frac{\vt}{\tvt} + \frac{\tvt}{\vt} - 2 \right) +  \tvt \Big[  s(\vr_b, \vt) -
\widetilde{s}  + \frac{\widetilde{e}}{\tvt} - \frac{e(\vr_b, \vt)}{\vt} \Big] \vr_b \vuB \cdot \vc{n}
- \widetilde{e} \left( \frac{\vt}{\tvt} - 1 \right) \vr_b \vuB \cdot \vc{n}
\right)
\D \sigma_x\\&\aleq \int_{\Gamma_{\rm in}} |\vt - \tvt|^2 \ \D \sigma_x.
\end{split}
\end{equation}
In view of \eqref{R9} and since the integrand is a sublinear function of $\vt$ for $\vt \to \infty$, it is enough to 
have
\[
\begin{split}
&\left( F_{i,b} \left( \frac{\vt}{\tvt} + \frac{\tvt}{\vt} - 2 \right) +  \tvt \Big[  s(\vr_b, \vt) -
\widetilde{s}  + \frac{\widetilde{e}}{\tvt} - \frac{e(\vr_b, \vt)}{\vt} \Big] \vr_b \vuB \cdot \vc{n}
- \widetilde{e} \left( \frac{\vt}{\tvt} - 1 \right) \vr_b \vuB \cdot \vc{n}
\right) \\ &\approx \tvt \left( \frac{F_{i,b}}{\vt} - \frac{e(\vr_b, \vt)}{\vt} \vr_b \vuB \cdot \vc{n} \right) 
\to - \infty \ \mbox{as}\ \vt \to 0.
\end{split}
\] 
{Indeed, in view of hypotheses \eqref{ws2}, \eqref{ws6},}
\[
\frac{ e(\vr_b, \vt) }{\vt} = \frac{1}{\vt} \frac{3}{2} \vr_b^{2/3} \frac{P \left( \frac{\vr_b} {\vt^{3/2}} \right)} 
{\left( \frac{\vr_b} {\vt^{3/2}} \right)^{5/3} }+ 
\frac{a}{\vr_b} \vt^3 \approx \frac{1}{\vt} \frac{3}{2} \vr^{2/3}_b p_\infty \ \mbox{as}\ \vt \to 0
\]
Moreover, by the same token,
\[
|s(\vr_b,\vt)| \aleq (1 - \log(\vt) )\ \mbox{as $\vt\to 0$}.
\]
Consequently, the desired estimate \eqref{ws13} holds as soon as 
\begin{equation} \label{ws14}
F_{i,b} - \frac{3}{2} p_\infty \vr_b^{5/3} \vuB \cdot \vc{n} < 0 \ \mbox{on} \ \Gamma_{\rm in}. 
\end{equation}
If \eqref{ws14} holds, the boundary integral in \eqref{R7} can be controlled via \eqref{ws11}, and we are in the situation 
treated in \cite[Section 6.1, Theorem 6.2]{EF101}. More precisely, a Gronwall type argument can be used to absorb all terms in 
${\rm Er}_2$ in \eqref{R7} to obtain the following result. 

\begin{Theorem}[Unconditional weak--strong uniqueness, I] \label{RT2}

Let the thermodynamic functions $p$, $e$, and $s$ satisfy the hypotheses \eqref{ws1}--\eqref{ws7}, 
where, in addition, 
\begin{equation} \label{ws3bis}
s(\vr, \vt) = \mathcal{S} \left( \frac{\vr} {\vt^{3/2}} \right) + \frac{4a}{3} \frac{\vt^3}{\vr},\ \mathcal{S}(Z) \to 0 \ \mbox{as}\ 
Z \to \infty.
\end{equation}
Let the transport 
coefficients $\mu$, $\eta$, $\kappa$ be continuously differentiable functions of $\vt \in [0, \infty)$, 
satisfying the hypotheses \eqref{ws8}--\eqref{ws10}.
Let the flux $F_{i,b}$ prescribed on the inflow boundary satisfy
\begin{equation} \label{ws14bis}
\sup_{x \in \Gamma_{\rm in} } \left( \frac{F_{i,b}}{|\vuB \cdot \vc{n}|}(x) + \frac{3}{2} p_\infty \vr_b^{5/3}(x) \right) < 0. 
\end{equation}
Let $(\vr, \vu, \vt)$ be a weak solution
of the Navier--Stokes system \eqref{i1}--\eqref{i4}, with the boundary conditions \eqref{I1}--\eqref{i10},  
in the sense specified in Definition \ref{wD1}.
Suppose that the same problem (with the same initial and boundary data) admits a strong solution $(\tvr, \tvu, \tvt)$ 
in the class 
\[
\tvr, \ \tvu, \ \tvt \in 
C^1([0,T] \times \Ov{\Omega}), \ \partial^2_x \tvu, \ \partial^2_x  \tvt \in 
C([0,T] \times \Ov{\Omega}). 
\]

Then 
\[
\vr = \tvr,\ \vu = \tvu, \ \vt = \tvt \ \mbox{in}\ [0,T] \times \Ov{\Omega}.
\]

\end{Theorem}

Hypothesis \eqref{ws14bis} may seem rather awkward, however, it can be interpreted as negativity of the ``heat flux'' on 
$\Gamma_{\rm in}$. To see this, consider an ``iconic'' example of internal energy satisfying \eqref{ws2}, namely, 
\begin{equation} \label{icon}
\vr e(\vr, \vt) = \underbrace{\frac{3}{2} \vr \vt}_{\rm molecular \ energy} + \underbrace{\frac{3}{2} p_\infty \vr^{5/3}}_{
\rm electron \ energy} + \underbrace{a \vt^4}_{\rm radiation (photon) energy}.
\end{equation}
Writing 
\[
F_{i,b} = \underbrace{\frac{3}{2} p_\infty \vr_b^{5/3} \vuB \cdot \vc{n}}_{\rm cold \ flux} + \underbrace{F^\tau_{i,b} \vu_b \cdot \vc{n}}_{\rm heat \ flux} \ \mbox{on}\ \Gamma_{\rm in}
\]
we can check that \eqref{ws14} holds as soon as $\inf_{\Gamma_{\rm in}} F^\tau_{i,b} > 0$.

{Although physically relevant, the satisfaction of the Third law may seem restrictive, in particular, 
this assumption is violated by the state equation \eqref{icon}. 
We show that the conclusion of Theorem \ref{RT2} remains valid in the general case under slightly more restrictive assumption} 
\begin{equation} \label{lambda}
\frac{1}{2} \leq \Lambda \leq 1 
\end{equation}
{in \eqref{ws8}. Note that the range \eqref{lambda} is still realistic for gases, see e.g. Becker \cite{BE}. As a matter of fact, 
the hypothesis \eqref{dod4} was not explicitly used in the proof of Theorem \ref{RT2}, it is necessary to control the error 
${\rm Er}_2$ exactly as in \cite[Section 6.1]{EF101}. A short inspection of the proof in \cite[Section 6.1]{EF101} reveals the most
problematic term in ${\rm Er}_2$, namely}  
\[
\intO{ \vr (\widetilde{s} - s)(\vu - \tvu) \cdot \Grad \tvt },
\]
{or, more precisely, its ``residual'' component}
\[
\int_{\mathcal{R} } \vr s (\vr, \vt) (\vu - \tvu) \cdot \Grad \tvt \dx, 
\]
{where}
\[
\mathcal{R} = \Big[ (0,T) \times \Omega \Big] \setminus \left\{ (t,x) \ \Big|\ \frac{1}{2} \inf \tvr \leq \vr(t,x) 
\leq 2 \sup \tvr,\ \frac{1}{2} \inf \tvt \leq \vt(t,x) 
\leq 2 \sup \tvt \right\},
\]
and where, by virtue of (\ref{ws1}--\ref{ws7}), 
\[
\|\vr 1_{\mathcal{ R}}\|_{L^{5/3}(\Omega)}^{5/3}\aleq E \left( \vr, \vu, S \Big|\ \tvr, \tvu, \widetilde{S} \right)
\]

The Gronwall argument used in \cite[Section 6.1]{EF101} applies as soon as we have the following bound:
\begin{equation} \label{estim}
\left| \int_{\mathcal{R} } \vr s (\vr, \vt) (\vu - \tvu) \cdot \Grad \tvt \dx \right| \leq 
\ep \left\| \vu - \tvu \right\|^2_{W^{1, \alpha}_0(\Omega; R^d)} + c(\ep) \intO{ \chi (t) E \left( \vr, \vu, S \Big|\ \tvr, \tvu, { \widetilde{S}} \right) } 
\end{equation}
for any $\ep > 0$, where
\[
\alpha = \frac{8}{5 - \Lambda}, \ \mbox{and} \ \chi \in L^1(0,T).
\]
In view of the arguments of \cite[Section 6.1]{EF101}, the verification of \eqref{estim} amounts to showing   
\begin{equation} \label{estima}
\int_{\mathcal{R} } \vr |\log(\vt)| |\vu - \tvu|  \dx  \leq 
\ep \left\| \vu - \tvu \right\|^2_{W^{1, \alpha}_0(\Omega; R^d)} + c(\ep) \chi (t)  \intO{  E \left( \vr, \vu, S \Big|\ \tvr, \tvu,  \widetilde{S} \right) }.
\end{equation}
Here, we address the problem for $d=3$, the result can be slightly improved for $d=2$. 
If $\Lambda \geq \frac{1}{2}$, we have $\alpha \geq \frac{16}{9}$, and 
in view of the standard Sobolev embedding theorem,
\[
W^{1, \alpha}_0(\Omega; R^3) \subset L^r(\Omega; R^3) \ \mbox{for}\ 1 \leq r \leq \frac{48}{11}.
\]
Consequently, by H\" older's inequality, 
\[
\begin{split}
\int_{\mathcal{R} } \vr |\log(\vt)| |\vu - \tvu|  \dx &\aleq 
\| 1_{\mathcal{R}} \vr \|_{L^{\frac{5}{3}}(\Omega)} \| \log(\vt) \|_{L^6(\Omega)} \| \vu - \tvu \|_{L^r(\Omega; R^3)}\\
&\leq \ep \left\| \vu - \tvu \right\|^2_{W^{1, \alpha}_0(\Omega; R^d)} + 
c(\ep)  \| \log(\vt) \|_{L^6(\Omega)}^2 \| 1_{\mathcal{R}} \vr \|^2_{L^{\frac{5}{3}}(\Omega)} \\
&\leq \ep \left\| \vu - \tvu \right\|^2_{W^{1, \alpha}_0(\Omega; R^d)} + 
c(\ep)  \| \log(\vt) \|_{L^6(\Omega)}^2  \left( \intO{ E \left( \vr, \vu, S \Big|\ \tvr, \tvu, {\widetilde{S}} \right) } \right)^{\frac{6}{5}}. 
\end{split}
\]
As $\vt$ is a weak solution, we have
\[
\log(\vt) \in L^2(0,T; W^{1,2}(\Omega)), \ \mbox{whence we may consider}\ 
\chi =  \| \log(\vt) \|_{L^6(\Omega)}^2 \in L^1(0,T).
\]

{We have shown the following result.} 
\begin{Theorem}[Unconditional weak--strong uniqueness, II] \label{RT3}

Let the thermodynamic functions $p$, $e$, and $s$ satisfy the hypotheses \eqref{ws1}--\eqref{ws7}. 
Let the transport 
coefficients $\mu$, $\eta$, $\kappa$ be continuously differentiable functions of $\vt \in [0, \infty)$, 
satisfying the hypotheses \eqref{ws8}--\eqref{ws10}, with 
\[
\Lambda \in \left[ \frac{1}{2}, 1 \right].
\]
Let the flux $F_{i,b}$ prescribed on the inflow boundary satisfy
\[
\sup_{x \in \Gamma_{\rm in} } \left( \frac{F_{i,b}}{|\vuB \cdot \vc{n}|}(x) + \frac{3}{2} p_\infty \vr_b^{5/3}(x) \right) < 0. 
\]
Let $(\vr, \vu, \vt)$ be a weak solution
of the Navier--Stokes system \eqref{i1}--\eqref{i4}, with the boundary conditions \eqref{I1}--\eqref{i10},  
in the sense specified in Definition \ref{wD1}.
Suppose that the same problem (with the same initial and boundary data) admits a strong solution $(\tvr, \tvu, \tvt)$ 
in the class 
\[
\tvr, \ \tvu, \ \tvt \in 
C^1([0,T] \times \Ov{\Omega}), \ \partial^2_x \tvu, \ \partial^2_x  \tvt \in 
C([0,T] \times \Ov{\Omega}). 
\]

Then 
\[
\vr = \tvr,\ \vu = \tvu, \ \vt = \tvt \ \mbox{in}\ [0,T] \times \Ov{\Omega}.
\]

\end{Theorem}

\section{Existence theory}
\label{E}

Our ultimate goal is to show existence of global--in--time weak solutions in the sense of Definition \ref{wD1}. To this end, we restrict ourselves to the thermodynamic functions $p$, $e$, and $s$, and the transport coefficients $\mu$, $\eta$, and $\kappa$ satisfying the constitutive relations \eqref{ws1}--\eqref{ws10} introduced in Section \ref{UWU}. The approximation scheme is similar to 
\cite[Chapter 3]{FeNo6A}, with the necessary modifications to accommodate the boundary conditions. In comparison with  
\cite[Chapter 3]{FeNo6A}, there are two main difficulties to be handled:
\begin{itemize}
\item Compactness of the boundary integrals with respect to the available {\it a priori} bounds. 
\item The fact that the approximate density does not satisfy the boundary condition \eqref{i8} at the first level of approximation. 

\end{itemize}

The existence result reads as follows:

\begin{Theorem}[Global--in--time existence] \label{tE1}
Let $\Omega \subset R^d$ be a bounded domain with smooth boundary. Let the thermodynamic functions $p$, $e$, and $s$ satisfy the hypotheses \eqref{ws1}--\eqref{ws7}, and let the transport coefficients be continuously differentiable functions of $\vt \in [0, \infty)$, 
satisfying the hypotheses \eqref{ws8}--\eqref{ws10}. Let the data $\vc{g}$, $\vuB$, $\vr_b$, and $F_{i,b}$ be smooth fucntions of 
$x \in \Ov{\Omega}$, satisfying 
\begin{equation} \label{E1}
\vr_b > 0 \ \mbox{on}\ \Gamma_{\rm in},\  
\sup_{x \in \Gamma_{\rm in} } \left( \frac{F_{i,b}}{|\vuB \cdot \vc{n}|}(x) + \frac{3}{2} p_\infty \vr_b^{5/3}(x) \right) < 0.
\end{equation}
Let the initial data $(\vr_0, (\vr \vu)_0, (\vr s)_0)$ be given such that 
\begin{equation} \label{E2}
\begin{split}
\vr(0, \cdot) &= \vr_0, \ \vr_0 \in L^{5/3}(\Omega), \ \vr_0 \geq 0,\ \mbox{a.a. in}\ \Omega; \\ 
(\vr \vu)(0, \cdot) &= (\vr \vu)_0,\ \intO{ \frac{|(\vr \vu)_0 |^2 }{\vr_0} } < \infty;\\
(\vr s)(0, \cdot) \equiv S(0, \cdot) &= (\vr s)_0 = \vr_0 s(\vr_0, \vt_0) \in L^1(\Omega),\ \vt_0 > 0 \ \mbox{a.a. in}\ \Omega.
\end{split}
\end{equation}

Then for any $T > 0$ the Navier--Stokes system \eqref{i1}--\eqref{i4}, with the boundary conditions \eqref{I1}--\eqref{i10}, 
admits a weak solution $(\vr, \vu, \vt)$ in $(0,T) \times \Omega$ in the sense of Definition \ref{wD1}. 

\end{Theorem}

The rest of the paper is devoted to the proof of Theorem \ref{tE1}.

\subsection{Approximation scheme}
\label{SE1}

Similarly to \cite[Chapter 3]{FeNo6A}, we introduce a multilevel approximation scheme to construct the weak solution. 

\subsubsection{Equation of continuity}

As in \cite{ChJiNo}, the equation of continuity is replaced by a standard parabolic regularization:
\begin{equation} \label{E3}
\begin{split}
\partial_t \vr + \Div (\vr \vu ) &= \ep \Del \vr \ \mbox{in}\ (0,T) \times \Omega,\ \ep > 0,\\
\ep \Grad \vr \cdot \vc{n} + (\vr_b - \vr) [\vuB \cdot \vc{n}]^- &= 0 \ \mbox{in}\ [0,T] \times \partial \Omega,\\
\vr(0, \cdot) &= \vr_{0, \delta},
\end{split}
\end{equation}
where we have denoted 
\[
[\vuB \cdot \vc{n}]^- = \min \left\{ 0 , \vuB \cdot \vc{n} \right\} = \left\{ \begin{array}{l} \vuB \cdot \vc{n}
\ \mbox{on}\ \Gamma_{\rm in} \\ 0 \ \mbox{otherwise} \end{array}               \right.
\]
Here $\ep > 0$, $\delta > 0$ are two parameters, $\vr_{0, \delta}$ being a suitable regularization of $\vr_0$,\ 
\[
\vr_{0, \delta} \in C^3(\Ov{\Omega}),\ 
\vr_{0, \delta} > 0 \ \mbox{in}\ \Ov{\Omega},\ \ep \Grad \vr_{0, \delta} \cdot \vc{n} + (\vr_b - \vr_{0, \delta})[ 
\vuB \cdot \vc{n}]^- = 0 \ \mbox{in}\ \partial \Omega.
\]

\subsubsection{Momentum equation}

The approximate velocities are determined via a Faedo--Galerkin approximation. To this end, consider  
\[
X_n = {\rm span} \left\{ \vc{w}_i\ \Big|\ \vc{w}_i \in \DC(\Omega; R^d),\ i = 1,\dots, n \right\}
\]
where $\vc{w}_i$ are orthonormal with respect to the standard scalar product in $L^2$. Let 
$\Pi_n: L^2 \to X_n$ be the associated orthogonal projection. 

We look for 
\[
\vu = \vc{v} + \vuB, \ \vc{v} \in C([0,T]; X_n),
\]
where 
\begin{equation} \label{E4}
\begin{split}
\left[ \intO{ \vr \vu \cdot \bfphi } \right]_{t=0}^{t = \tau} &=
\int_0^\tau \intO{ \Big[ \vr \vu \cdot \partial_t \bfphi + \vr \vu \otimes \vu : \Grad \bfphi
+ p_\delta \Div \bfphi - \mathbb{S}_\delta   : \Grad \bfphi \Big] }\\
&- \ep \int_0^\tau \intO{ \Grad \vr \cdot \Grad \vu \cdot \bfphi } \dt + \int_0^\tau \intO{ \vr \vc{g} \cdot \bfphi } \dt
\end{split}
\end{equation}
for any $\bfphi \in C^1([0,T]; X_n)$, with the initial condition $\Pi_n (\vr \vu)_0$. Here we have introduced 
\[
\begin{split}
p_\delta &= p + \delta \left(\vr^\Gamma + \vr^2 \right),\\ 
\mathbb{S}_\delta (\vt, \Grad \vu) &= 
(\mu(\vt) + \delta \vt) \left( \Grad \vu + \Grad^t \vu - \frac{2}{d} \Div \vu \mathbb{I} \right) + 
\eta (\vt) \ \Div \vu \mathbb{I}.  
\end{split}
\]

\subsubsection{Internal energy balance}

To keep the approximate scheme consistent with the energy estimates, we consider a modified internal energy balance:
\begin{equation} \label{E5}
\partial_t (\vr e_\delta) + \Div (\vr e_\delta \vu) + \Div \vc{q}_\delta =
\mathbb{S}_\delta  : \Grad \vu - p \Div \vu + \ep \delta \left( \Gamma \vr^{\Gamma - 2} + 2 \right) |\Grad \vr|^2 
+ \delta \frac{1}{\vt^2} - \ep \vt^5, 
\end{equation}
with the Robin boundary conditions 
\begin{equation} \label{E6}
\vc{q}_\delta  \cdot \vc{n} +
\vr e_\delta [\vuB \cdot \vc{n}]^- = F_{i,b}, \ F_{i,b} = 0 \ \mbox{on} \ \Gamma_{\rm wall} \cup \Gamma_{\rm out}.
\end{equation}
Here, 
\[
e_\delta = e + \delta \vt, \ \vc{q}_\delta = \vc{q} - \delta \left( \vt^\Gamma + \frac{1}{\vt} \right) \Grad \vt.
\]
The initial conditions are determined through $\vr_0$ and $\vt_0$, where 
\begin{equation} \label{E6a}
\vt_{0, \delta} \in W^{1,2} \cap L^\infty (\Omega), \ {\rm ess} \inf_{x \in \Omega} \vt_{0, \delta}(x) > 0.
\end{equation}

\subsubsection{Entropy inequality}

The above approximation scheme is exactly the same as in \cite[Chapter 3]{FeNo6A}, modulo the boundary conditions. Assuming 
there is a strong solution of the internal energy balance such that 
\[
\inf_{(0,T) \times \Omega} \vt > 0,
\] 
we derive a \emph{weak formulation} of the entropy inequality. 

First, rewrite \eqref{E5} with the help of \eqref{E3} as 
\[
\begin{split}
\vr \frac{1}{\vt} \partial_t
&e_\delta + \frac{1}{\vt}\vr \vu \cdot \Grad e_\delta \vu  + \ep e_\delta \frac{1}{\vt} \Del \vr + 
\Div \left( \frac{\vc{q}_\delta}{\vt} \right) \\&=
\frac{1}{\vt} \left( \mathbb{S}_\delta  : \Grad \vu - \frac{\vc{q}_\delta \cdot \Grad \vt}{\vt} + \frac{\delta}{\vt^2} \right)  - \frac{p}{\vt} \Div \vu + \ep \delta \frac{1}{\vt} \left( \Gamma \vr^{\Gamma - 2} + 2 \right) |\Grad \vr|^2 
 - \ep \vt^4. 
\end{split}
\]
Consequently, using Gibbs' relation, we may go back to the entropy formulation:
\begin{equation} \label{E7}
\begin{split}
\partial_t  &(\vr s_\delta) + \Div (\vr s_\delta \vu )+  
\Div \left( \frac{\vc{q}_\delta}{\vt} \right) \\&=
\frac{1}{\vt} \left( \mathbb{S}_\delta  : \Grad \vu - \frac{\vc{q}_\delta \cdot \Grad \vt}{\vt} + \frac{\delta}{\vt^2} \right)  + \ep \delta \frac{1}{\vt} \left( \Gamma \vr^{\Gamma - 2} + 2 \right) |\Grad \vr|^2 
 - \ep \vt^4 \\ 
&- \ep \Del \vr \left( \frac{e_\delta}{\vt} - s_\delta + \frac{p}{\vr \vt} \right), 
\end{split}
\end{equation}
where we have denoted 
\[
s_\delta = s + \delta \log(\vt).
\]

Finally, multiplying \eqref{E7} by $\varphi \in C^1([0,T] \times \Ov{\Omega})$, $\varphi \geq 0$, we get 
\[
\begin{split}
& \left[ \intO{ \vr s_\delta \varphi } \right]_{t=0}^{t = \tau} - \int_0^\tau \intO{ \vr s_\delta \partial_t \varphi } \dt 
- \int_0^\tau \intO{ \vr s_\delta \vu \cdot \Grad \varphi } \dt 
 + 
\int_0^\tau \int_{\Gamma_{\rm out}} \varphi \vr s_\delta \vuB \cdot \vc{n} \ \D \sigma_x \dt \\& -  
\int_0^\tau \intO{ \frac{\vc{q}_\delta}{\vt} \cdot \Grad \varphi } \dt \\&=
\int_0^\tau \intO{ \frac{\varphi}{\vt} \left( \mathbb{S}_\delta  : \Grad \vu - \frac{\vc{q}_\delta \cdot \Grad \vt}{\vt} + \frac{\delta}{\vt^2} \right)} \dt  + \ep \int_0^\tau \intO{ \delta \frac{\varphi}{\vt} \Big( \left( \Gamma \vr^{\Gamma - 2} + 2 \right) |\Grad \vr|^2  - \vt^4 \Big) } \dt \\
&+ \ep \int_0^\tau \intO{ \Grad \vr  \cdot \Grad \left[ \varphi  \left( \frac{e_\delta}{\vt} - s_\delta + \frac{p}{\vr \vt} \right)
\right] } \dt \\ 
&- \ep \int_0^\tau \int_{\Gamma_{\rm in}} \varphi \Grad \vr \cdot \vc{n} \left( \frac{e_\delta}{\vt} - s_\delta + \frac{p}{\vr \vt} \right)
\ \D \sigma_x \dt - \int_0^\tau \int_{\Gamma_{\rm in}} 
\varphi \left( \frac{\vc{q}_\delta \cdot \vc{n}}{\vt}
+ \vr s_\delta \vuB \cdot \vc{n} \right) \D \sigma_x \dt 
\end{split}
\]

Now, we use the boundary conditions obtaining
\[
\begin{split}
- \ep &\int_0^\tau \int_{\Gamma_{\rm in}} \varphi \Grad \vr \cdot \vc{n} \left( \frac{e_\delta}{\vt} - s_\delta + \frac{p}{\vr \vt} \right)
\ \D \sigma_x \dt - \int_0^\tau \int_{\Gamma_{\rm in}} 
\varphi \left( \frac{\vc{q}_\delta \cdot \vc{n}}{\vt}
+ \vr s_\delta \vuB \cdot \vc{n} \right) \D \sigma_x \dt 
\\
&=
\int_0^\tau \int_{\Gamma_{\rm in}} \varphi \left[ \vr_b \left(\frac{e_\delta}{\vt} - s_\delta + \frac{p}{\vr \vt} \right) 
\vuB \cdot \vc{n} - \frac{p}{\vt} \vuB \cdot \vc{n}  - \frac{1}{\vt} F_{i,b} \right] \D \sigma_x \dt.
\end{split}
\] 
Unfortunately, the boundary integrals depend on $\vr$ on the inflow part $\Gamma_{\rm in}$ whereas, at this stage of approximation, 
$\vr \ne \vr_b$ on $\Gamma_{\rm in}$, in general. To remedy the problem, we make use
of the thermodynamic stability \eqref{i5}, \eqref{i5b} to deduce
\[
\begin{split}
E_{\rm int}(\vr_b, S_b) &\geq E_{\rm int}(\vr, S) + \left( e - \vartheta s + \frac{p}{\vr} \right) (\vr_b - \vr) +
\vt (S_b - S) \\
&= \left( e - \vartheta s + \frac{p}{\vr} \right) \vr_b + \vr e - \left( e - \vartheta s + \frac{p}{\vr} \right) \vr
+ \vt S_b - \vr \vt s\\
&= \left( e - \vartheta s + \frac{p}{\vr} \right) \vr_b - p + \vt S_b, \ \mbox{with}\ S_b = \vr_b s(\vr_b, \vt).
\end{split}
\]
Thus we may infer that
\[
\frac{1}{\vt}\left( e - \vartheta s + \frac{p}{\vr} \right) \vr_b - \frac{1}{\vt} p \leq \frac{1}{\vt} \vr_b e(\vr_b, \vt) - \vr_b s(\vr_b, \vt)
\]
Seeing that $\vuB \cdot \vc{n} < 0$ on $\Gamma_{\rm in}$ we may
write down the final form of the approximate entropy inequality: 
\begin{equation} \label{E8}
\begin{split}
& \left[ \intO{ \vr s_\delta \varphi } \right]_{t=0}^{t = \tau} - \int_0^\tau \intO{ \vr s_\delta \partial_t \varphi } \dt 
- \int_0^\tau \intO{ \vr s_\delta \vu \cdot \Grad \varphi } \dt 
 + 
\int_0^\tau \int_{\Gamma_{\rm out}} \varphi \vr s_\delta \vuB \cdot \vc{n} \ \D \sigma_x \dt \\& -  
\int_0^\tau \intO{ \frac{\vc{q}_\delta}{\vt} \cdot \Grad \varphi } \dt \\& \geq 
\int_0^\tau \intO{ \frac{\varphi}{\vt} \left( \mathbb{S}_\delta  : \Grad \vu - \frac{\vc{q}_\delta \cdot \Grad \vt}{\vt} + \frac{\delta}{\vt^2} \right)} \dt  + \ep \int_0^\tau \intO{ \delta \frac{\varphi}{\vt} \Big( \left( \Gamma \vr^{\Gamma - 2} + 2 \right) |\Grad \vr|^2  - \vt^4 \Big) } \dt \\
&+ \ep \int_0^\tau \intO{ \Grad \vr  \cdot \Grad \left[ \varphi  \left( \frac{e_\delta}{\vt} - s_\delta + \frac{p}{\vr \vt} \right)
\right] } \dt \\
&+ \int_0^\tau \int_{\Gamma_{\rm in}} \varphi \left[ \delta \vr_b (1 - \log(\vt) ) 
\vuB \cdot \vc{n} - \frac{1}{\vt} F_{i,b} + 
\left( \frac{e(\vr_b, \vt)}{\vt} - s(\vr_b, \vt) \right) \vr_b \vuB \cdot \vc{n}  \right] \D \sigma_x \dt  
\end{split}
\end{equation}
for any $\varphi \in C^1([0,T] \times \Ov{\Omega})$, $\varphi \geq 0$.

\subsection{A priori bounds}

Assuming for a moment solvability of the approximate problem \eqref{E1}--\eqref{E6a} we focus on available {\it a priori} bounds. 
Needless to say they mimick their counterparts for the limit system. 

\subsubsection{Approximate total energy balance}

Consider $\bfphi = \psi \vc{v} = \psi (\vu - \vuB)$, $\psi = \psi(t)$, $\psi \in C^1[0,T]$ as a test function in the approximate momentum balance \eqref{E4}. After a straightforward 
manipulation and with the help of the equation of continuity, we deduce 
\begin{equation} \label{E10}
\begin{split}
&\left[ \intO{\left[ \frac{1}{2} \vr |\vu - \vuB|^2 + \delta \left( \frac{1}{\Gamma - 1} \vr^\Gamma + \vr^2 \right) \right] \psi } \right]_{t = 0}^{ t = \tau}\\ &- \int_0^\tau \partial_t \psi \intO{ \left[ \frac{1}{2} \vr |\vu - \vuB|^2 + \delta \left( \frac{1}{\Gamma - 1} \vr^\Gamma + \vr^2 \right) \right]} \dt + 
\int_0^\tau \psi \intO{ \mathbb{S}_\delta (\vt, \Grad \vu) : \Grad \vu } \dt \\ 
&+\int_0^\tau \psi \int_{\Gamma_{\rm out}} \delta \left( \frac{1}{\Gamma - 1} \vr^\Gamma + \vr^2 \right) \vuB \cdot \vc{n} \ \D \sigma_x \dt\\
&- \delta \int_0^\tau \psi \int_{\Gamma_{\rm in}} \left[ \frac{1}{\Gamma - 1} \vr_b^\Gamma - \frac{\Gamma}{\Gamma - 1}\vr^{\Gamma - 1} (\vr_b - \vr) - \frac{1}{\Gamma - 1}\vr^\Gamma \right] \vuB \cdot \vc{n} 
\ {\rm d}\sigma_x \ \dt\\
&- \delta \int_0^\tau \psi \int_{\Gamma_{\rm in}} (\vr - \vr_b )^2 \vu_B \cdot \vc{n} 
\ {\rm d}\sigma_x \ \dt + \ep \delta \int_0^\tau \psi \intO{ \left( \Gamma \vr^{\Gamma - 2} + 2 \right)  |\Grad \vr|^2 } \dt
\end{split}
\end{equation}
\[
\begin{split}
=
&- 
\int_0^\tau \psi \intO{ \left[ \vr \vu \otimes \vu + p_\delta (\vr) \mathbb{I} \right]  :  \Grad \vu_B } \dt + \frac{1}{2} \int_0^\tau 
\psi {\intO{ {\vr} \vu  \cdot\Grad |\vu_B|^2  } }
\dt\\ &+ \int_0^\tau \psi \intO{ p \Div \vu } \dt 
\\ &+ \int_0^\tau \psi \intO{ \mathbb{S}_\delta(\vt,  \Grad \vu) : \Grad \vu_B } \dt + 
\int_0^\tau \psi \intO{ \vr \vc{g} \cdot (\vu - \vuB) } \dt  \\&- \delta \int_0^\tau \psi \int_{\Gamma_{\rm in}} \frac{1}{\Gamma - 1}\vr_b^\Gamma  \vuB \cdot \vc{n} \ \D \sigma_x \dt + \ep \int_0^\tau \psi \intO{ \Grad \vr \cdot \Grad (\vu - \vuB) \cdot \vuB } \dt  
\end{split}
\]
Equality \eqref{E10} added to the approximate internal energy balance \eqref{E5}, 
\eqref{E6} gives rise to the {\em approximate total energy balance} 
\[
\begin{split}
&\left[ \intO{\left[ \frac{1}{2} \vr |\vu - \vuB|^2 + \delta \left( \frac{1}{\Gamma - 1} \vr^\Gamma + \vr^2 \right) + 
\vr e_\delta  \right] \psi } \right]_{t = 0}^{ t = \tau}\\ 
&- \int_0^\tau \partial_t \psi \intO{ \left[ \frac{1}{2} \vr |\vu - \vuB|^2 + \delta \left( \frac{1}{\Gamma - 1} \vr^\Gamma + \vr^2 \right) + 
\vr e_\delta  \right] } \dt  
\\ 
&+\int_0^\tau \psi \int_{\Gamma_{\rm out}} \vr e_\delta \vuB \cdot \vc{n} \D \sigma_x \dt   + \delta \int_0^\tau \int_{\Gamma_{\rm out}}  \left( \frac{1}{\Gamma - 1} \vr^\Gamma + \vr^2 \right) \vuB \cdot \vc{n} \ \D \sigma_x \dt\\
&+\int_0^\tau \psi \int_{\Gamma_{\rm in}} F_{i,b} \D \sigma_x \dt- \delta \int_0^\tau \psi \int_{\Gamma_{\rm in}} \left[ \frac{1}{\Gamma - 1} \vr_b^\Gamma - \frac{\Gamma}{\Gamma - 1}\vr^{\Gamma - 1} (\vr_b - \vr) - \frac{1}{\Gamma - 1}\vr^\Gamma \right] \vuB \cdot \vc{n} 
\ {\rm d}\sigma_x \ \dt  \\
&- \delta \int_0^\tau \psi \int_{\Gamma_{\rm in}} (\vr - \vr_b )^2 \vu_B \cdot \vc{n} 
\ {\rm d}\sigma_x \ \dt 
\end{split}
\]
\begin{equation} \label{E11}
\begin{split}
&= - 
\int_0^\tau \psi \intO{ \left[ \vr \vu \otimes \vu + p_\delta  \mathbb{I} \right]  :  \Grad \vu_B } \dt + \frac{1}{2} \int_0^\tau 
\psi {\intO{ {\vr} \vu  \cdot\Grad |\vu_B|^2  } }
\dt  
\\ &+ \int_0^\tau \psi \intO{ \mathbb{S}_\delta(\vt,  \Grad \vu) : \Grad \vuB } \dt + 
\int_0^\tau \psi \intO{ \vr \vc{g} \cdot (\vu - \vuB) } \dt\\ &+ \int_0^\tau \psi \intO{ \left( \delta \frac{1}{\vt^2} - 
\ep \vt^5 \right)} \dt \\&- \delta \int_0^\tau \psi \int_{\Gamma_{\rm in}} \frac{1}{\Gamma - 1}\vr_b^\Gamma  \vuB \cdot \vc{n} \ \D \sigma_x \dt  + \ep \int_0^\tau \psi \intO{ \Grad \vr \cdot \Grad (\vu - \vuB) \cdot \vuB } \dt 
\end{split}
\end{equation}
for any $\psi \in C^1[0,T]$.

\subsubsection{Approximate total entropy balance and uniform bounds}	

Consider $\varphi \Ov{\vt} > 0$, $\Ov{\vt}$ --a positive constant--  as a test function in the approximate entropy inequality \eqref{E8} and subtract the resulting expression from \eqref{E11} (with $\psi \equiv 1$) to obtain
\begin{equation} \label{E12}
\begin{split}
&\left[ \intO{\left[ \frac{1}{2} \vr |\vu - \vuB|^2 + \delta \left( \frac{1}{\Gamma - 1} \vr^\Gamma + \vr^2 \right) + 
\vr e_\delta - \Ov{\vt} \vr s_\delta \right] } \right]_{t = 0}^{ t = \tau} 
 \\ 
&+\int_0^\tau \int_{\Gamma_{\rm out}} \vr (e_\delta - \Ov{\vt} s_\delta) \vuB \cdot \vc{n} \D \sigma_x \dt   + \delta \int_0^\tau \int_{\Gamma_{\rm out}}  \left( \frac{1}{\Gamma - 1} \vr^\Gamma + \vr^2 \right) \vuB \cdot \vc{n} \ \D \sigma_x \dt\\
&+\int_0^\tau \int_{\Gamma_{\rm in}}  \left[ \left( 1 - \frac{\Ov{\vt}}{\vt} \right) F_{i,b} + 
\Ov{\vt} \left( \frac{e(\vr_b, \vt)}{\vt} - s(\vr_b, \vt) \right) \vr_b \vuB \cdot \vc{n}  \right] \D \sigma_x \dt 
\\
&- \delta \int_0^\tau \int_{\Gamma_{\rm in}} \left[ \frac{1}{\Gamma - 1} \vr_b^\Gamma - \frac{\Gamma}{\Gamma - 1}\vr^{\Gamma - 1} (\vr_b - \vr) - \frac{1}{\Gamma - 1}\vr^\Gamma \right] \vuB \cdot \vc{n} 
\ {\rm d}\sigma_x \ \dt  \\
&- \delta \int_0^\tau \int_{\Gamma_{\rm in}} (\vr - \vr_b )^2 \vuB \cdot \vc{n} 
\ {\rm d}\sigma_x \ \dt + \delta \int_0^\tau \int_{\Gamma_{\rm in}} \Ov{\vt} \vr_b (1 - \log(\vt)) \vuB \cdot \vc{n} \D \sigma_x \dt
\\
&+ \int_0^\tau \intO{ \frac{\Ov{\vt}}{\vt} \left( \mathbb{S}_\delta  : \Grad \vu - \frac{\vc{q}_\delta \cdot \Grad \vt}{\vt} + \frac{\delta}{\vt^2} \right)} \dt  + \ep \delta \int_0^\tau \intO{  \frac{\Ov{\vt}}{\vt}  \left( \Gamma \vr^{\Gamma - 2} + 2 \right) |\Grad \vr|^2  } \dt
\\
&+ \ep \int_0^\tau \intO{ \left( \vt^5 - \Ov{\vt} \vt^4 \right) } \dt - \delta \int_0^\tau \intO{ \frac{1}{\vt^2} }
\\
&\leq - 
\int_0^\tau \intO{ \left[ \vr \vu \otimes \vu + p_\delta (\vr) \mathbb{I} \right]  :  \Grad \vuB } \dt + \frac{1}{2} \int_0^\tau {\intO{ {\vr} \vu  \cdot\Grad |\vuB|^2  } }
\dt  
\\ &+ \int_0^\tau \intO{ \mathbb{S}_\delta(\vt,  \Grad \vu) : \Grad \vuB } \dt + 
\int_0^\tau \intO{ \vr \vc{g} \cdot (\vu - \vuB) } \dt  \ \\&- \delta \int_0^\tau \int_{\Gamma_{\rm in}} \frac{1}{\Gamma - 1}\vr_b^\Gamma  \vuB \cdot \vc{n} \ \D \sigma_x \dt 
-\ep \Ov{\vt} \int_0^\tau \intO{ \Grad \vr  \cdot \Grad \left[ \left( \frac{e_\delta}{\vt} - s_\delta + \frac{p}{\vr \vt} \right)
\right] } \dt \\ &+\ep \int_0^\tau \intO{ \Grad \vr \cdot \Grad (\vu - \vuB) \cdot \vuB } \dt 
\end{split}
\end{equation}

Now, the first observation is that the left--hand side of \eqref{E12} is bounded below by a constant that depends only on the 
data but is independent of $\ep$ and $\delta$. Indeed, similarly to Section \ref{UWU}, we may control the boundary integral 
\[
\int_{\Gamma_{\rm in}}  \left[ \left( 1 - \frac{\Ov{\vt}}{\vt} \right) F_{i,b} + 
\Ov{\vt} \left( \frac{e(\vr_b, \vt)}{\vt} - s(\vr_b, \vt) \right) \vr_b \vuB \cdot \vc{n}  \right] \D \sigma_x
\]
with the help of hypothesis \eqref{E1}:
\begin{equation} \label{E13}
\int_{\Gamma_{\rm in}} \left( \frac{1}{\vt} + \vt^3 |\vuB \cdot \vc{n}| \right) \D \sigma_x  \aleq 
\int_{\Gamma_{\rm in}}  \left[ \left( 1 - \frac{\Ov{\vt}}{\vt} \right) F_{i,b} + 
\Ov{\vt} \left( \frac{e(\vr_b, \vt)}{\vt} - s(\vr_b, \vt) \right) \vr_b \vuB \cdot \vc{n}  \right] \D \sigma_x + 1.
\end{equation}
The remaining boundary integrals are either non--negative or controllable by the quantity on the left-hand side in \eqref{E13}. 
Note that $\vuB \cdot \vc{n} < 0$ on $\Gamma_{\rm in}$ and $\vuB \cdot \vc{n} > 0$ on $\Gamma_{\rm out}$. 

The next observation is that all integrals on the right--hand side of \eqref{E12} can be ``absorbed'' in the left hand--side by means of a Gronwall argument. Indeed possibly the most difficult term is 
\[
\begin{split}
&\intO{ \mathbb{S}_\delta (\vt, \Grad \vu) : \Grad \vuB }\\ &= 
\intO{ (\mu(\vt) + \delta \vt ) \left( \Grad \vu + \Grad \vu^t - \frac{2}{d} \Div \vu \mathbb{I} \right) : \Grad \vuB }
+ \intO{ \eta (\vt) \Div \vu \Div \vuB }\\ 
&\leq \omega \intO{ \frac{1}{\vt} \mathbb{S}_\delta : \Grad \vu } + c_1(\omega, \vuB) \intO{ \vt^2 } 
\leq \omega \intO{ \frac{1}{\vt} \mathbb{S}_\delta : \Grad \vu } + c_2 (\omega, \vuB) \left(1 + \intO{ \vr e_\delta } \right)
\end{split}
\]
for any $\omega > 0$. 

As for the last integral in \eqref{E12}, it can be shown, exactly as in \cite[Chapter 3]{FeNo6A}, that
\begin{equation} \label{E14} 
\ep \Ov{\vt} \int_0^\tau \intO{ \Grad \vr  \cdot \Grad \left[ \left( \frac{e_\delta}{\vt} - s_\delta + \frac{p}{\vr \vt} \right)
\right] } \dt \to 0 \ \mbox{as}\ \ep \to 0
\end{equation}
for any fixed $\delta \to 0$.

Summarizing we have obtained the following {\it a priori} bounds, cf. \cite[Chapter 3]{FeNo6A}:\footnote{In what follows, we denote $a\aleq b$ if there exists $c>0$ independent of $n$, $\ep$, $\delta$ such that $a\le c b$.}
\begin{equation} \label{E15}
\begin{split}
{\rm ess} \sup_{t \in (0,T)} \intO{\left[ \frac{1}{2} \vr |\vu - \vuB|^2 + \delta \left( \frac{1}{\Gamma - 1} \vr^\Gamma + \vr^2 \right) + 
\vr e_\delta - \Ov{\vt} \vr s_\delta \right] } &\aleq 1,\\ 
\int_0^T \intO{ \frac{\Ov{\vt}}{\vt} \left( \mathbb{S}_\delta  : \Grad \vu - \frac{\vc{q}_\delta \cdot \Grad \vt}{\vt} \right)} \dt
&\aleq 1, \\ 
\int_0^T \int_{\Gamma_{\rm in}} \left( \frac{1}{\vt} + \vt^3 |\vuB \cdot \vc{n}| \right) \D \sigma_x \dt &\aleq 1,\\
\int_0^T \int_{\Gamma_{\rm out}} \vr (e_\delta - \Ov{\vt} s_\delta) \vuB \cdot \vc{n} \D \sigma_x \dt &\aleq 1,\\ 
\end{split}
\end{equation} 
and 
\begin{equation} \label{E16}
\begin{split}
\delta \int_0^T \intO{ \frac{1}{\vt^3} } \dt + \ep \int_0^T \intO{ \vt^5 } \dt &\aleq 1, \\ 
\delta \left( \int_0^\tau \int_{\Gamma_{\rm out}}  \left( \frac{1}{\Gamma - 1} \vr^\Gamma + \vr^2 \right) |\vuB \cdot \vc{n}| \ \D \sigma_x \dt + \int_0^T \int_{\Gamma_{\rm in}} (\vr - \vr_b)^2 |\vuB \cdot \vc{n}| \D \sigma_x \dt  \right) 
&\aleq 1, \\ 
\ep \delta \int_0^T \intO{ \frac{{1}}{\vt}  \left( \Gamma \vr^{\Gamma - 2} + 2 \right) |\Grad \vr|^2  } \dt &\aleq 1.
\end{split}
\end{equation}

\subsection{Solvability of the approximate problem}

Given $n \geq 1$, $\ep > 0$, and $\delta > 0$, the existence of solutions to the approximate system \eqref{E3}--\eqref{E6a}
was shown in \cite[Chapter 3]{FeNo6} in the case of energetically insulated system $\vuB \equiv 0$. The scheme of the proof 
is based on a fixed point argument:
\begin{enumerate}
\item Fix $[\vr, \vc{v}, \vt]$. 
\item Solve the approximate equation of continuity obtaining a new density $\vr = \vr[\vc{v}]$.
\item For given $\vc{v}, \vr[\vc{v}]$ solve the approximate internal energy equation \eqref{E5}--\eqref{E6a}  
$\vt[\vr, \vu]$
\item Find a new velocity $\vu$ by solving \eqref{E4} and use a fixed point argument. 
\end{enumerate}

In the present setting, the steps 1,2 have been performed in \cite{ChJiNo}. We therefore focus on solvability of the 
approximate internal energy equation \eqref{E5}--\eqref{E6a} for given (smooth) $\vr$ and $\vu$. This amounts to verifying 
the same set of {\it a priori} estimates as in \cite[Chapter 3, Section 3.4.2] {FeNo6A}. For given 
$\vr$, $\vu$, we consider the problem
\begin{equation} \label{eE5}
\partial_t (\vr e_\delta) + \Div (\vr e_\delta \vu) + \Div \vc{q}_\delta =
\mathbb{S}_\delta  : \Grad \vu - p \Div \vu + \ep \delta \left( \Gamma \vr^{\Gamma - 2} + 2 \right) |\Grad \vr|^2 
+ \delta \frac{1}{\vt^2} - \ep \vt^5, 
\end{equation}
\begin{equation} \label{eE6}
\vc{q}_\delta  \cdot \vc{n} +
\vr e_\delta [\vuB \cdot \vc{n}]^- = F_{i,b}, \ F_{i,b} = 0 \ \mbox{on} \ \Gamma_{\rm wall} \cup \Gamma_{\rm out},
\end{equation}
with
\[ 
\vc{q}_\delta = - \kappa (\vt) \Grad \vt - \delta \left( \vt^\Gamma + \frac{1}{\vt} \right) \Grad \vt.
\]
\subsubsection{Comparison principle}

We say that $\Ov{\vt}$ is a \emph{supersolution} of \eqref{eE5}, \eqref{eE6} if it satisfies 
\[
\partial_t (\vr e_\delta) + \Div (\vr e_\delta \vu) + \Div \vc{q}_\delta \geq 
\mathbb{S}_\delta  : \Grad \vu - p \Div \vu + \ep \delta \left( \Gamma \vr^{\Gamma - 2} + 2 \right) |\Grad \vr|^2 
+ \delta \frac{1}{\vt^2} - \ep \vt^5
\]
\[
\vc{q}_\delta  \cdot \vc{n} +
\vr e_\delta \vuB \cdot \vc{n} \leq F_{i,b} \ \mbox{on}\ \Gamma_{\rm in}, \ \vc{q}_\delta = 0 \ \mbox{on} \ \Gamma_{\rm wall} \cup \Gamma_{\rm out}.
\] 
Similarly, a \emph{subsolution} $\underline{\vt}$ satisfies
\[
\partial_t (\vr e_\delta) + \Div (\vr e_\delta \vu) + \Div \vc{q}_\delta \leq 
\mathbb{S}_\delta  : \Grad \vu - p \Div \vu + \ep \delta \left( \Gamma \vr^{\Gamma - 2} + 2 \right) |\Grad \vr|^2 
+ \delta \frac{1}{\vt^2} - \ep \vt^5
\]
\[
\vc{q}_\delta  \cdot \vc{n} +
\vr e_\delta \vuB \cdot \vc{n} \geq F_{i,b} \ \mbox{on}\ \Gamma_{\rm in}, \ \vc{q}_\delta = 0 \ \mbox{on} \ \Gamma_{\rm wall} \cup \Gamma_{\rm out}.
\] 

The \emph{comparison principle} asserts that if $\Ov{\vt}$ is a supersolution and $\underline{\vt}$ a subsolution, then
\begin{equation} \label{CP}
\underline{\vt}(0,  \cdot) \leq \Ov{\vt}(0, \cdot) \ \Rightarrow \ 
\underline{\vt}(t,  \cdot) \leq \Ov{\vt}(t, \cdot) \ \mbox{for all}\ t \geq 0.
\end{equation}

Our goal is establish the comparison principle for strong solutions of  the problem \eqref{eE5}, \eqref{eE6}. Following the proof in the case of homogeneous boundary conditions in \cite[Chapter 3, Lemma 3.2]{FeNo6A}, we consider the difference of the two inqualities multiplied 
by the expression 
\[
{\rm sgn}^+ \Big( \vr e_\delta (\vr, \underline{\vt}) - \vr e_\delta (\vr, \Ov{\vt}) \Big).    
\]
Specifically, we get 
\[
\begin{split} 
\left[ \partial_t \Big( \vr e_\delta (\vr, \underline{\vt}) - \vr e_\delta (\vr, \Ov{\vt}) \Big) + 
\vu \cdot \Grad \Big( \vr e_\delta (\vr, \underline{\vt}) - \vr e_\delta (\vr, \Ov{\vt}) \Big) \right] 
{\rm sgn}^+ \Big( \vr e_\delta (\vr, \underline{\vt}) - \vr e_\delta (\vr, \Ov{\vt}) \Big)\\
+ \Div \Big( \vc{q}_\delta(\underline{\vt}) -  \vc{q}_\delta(\Ov{\vt}) \Big){\rm sgn}^+ \Big( \vr e_\delta (\vr, \underline{\vt}) - \vr e_\delta (\vr, \Ov{\vt}) \Big) \leq \dots
\end{split}
\]
Thus, in comparison with \cite[Chapter 3, Lemma 3.2]{FeNo6A}, there is an extra term on the left-hand side of the above inequality
after integration, namely 
\[
\begin{split}
\int_{\partial \Omega} &{\rm sgn}^+ \Big( \vr e_\delta (\vr, \underline{\vt}) - \vr e_\delta (\vr, \Ov{\vt}) \Big) \left[ 
\Big( \vr e_\delta (\vr, \underline{\vt}) - \vr e_\delta (\vr, \Ov{\vt}) \Big) \vuB \cdot \vc{n} +
\Big( \vc{q}_\delta(\underline{\vt}) -  \vc{q}_\delta(\Ov{\vt}) \Big) \cdot \vc{n} \right] \D \sigma_x\\
& \geq \int_{\Gamma_{\rm in}} {\rm sgn}^+ \Big( \vr e_\delta (\vr, \underline{\vt}) - \vr e_\delta (\vr, \Ov{\vt}) \Big) \left[ 
\Big( \vr e_\delta (\vr, \underline{\vt}) - \vr e_\delta (\vr, \Ov{\vt}) \Big) \vuB \cdot \vc{n} +
\Big( \vc{q}_\delta(\underline{\vt}) -  \vc{q}_\delta(\Ov{\vt}) \Big) \cdot \vc{n} \right] \D \sigma_x\\ 
&\geq 0.
\end{split}
\]
Here, similarly to the proof in \cite[Chapter 3, Lemma 3.2]{FeNo6A}, we have written 
\[
\vc{q}_\delta (\vt) = - \Grad \mathcal{K}_\delta (\vt),\ \mathcal{K}_\delta (\vt) \equiv \int_1^\vt \left( \kappa(s) + \delta \left( s^\Gamma + \frac{1}{s} \right) \right) \D s,
\] 
and used the equality 
\[
{\rm sgn}^+ \Big( \vr e_\delta (\vr, \underline{\vt}) - \vr e_\delta (\vr, \Ov{\vt}) \Big) = 
{\rm sgn}^+ \Big( \mathcal{K}_\delta (\underline{\vt})  - \mathcal{K}_\delta (\Ov{\vt}) \Big).
\]
Accordingly, the proof of \eqref{CP} can be carried over exactly as in \cite[Chapter 3, Lemma 3.2]{FeNo6A}.

As a corollary, we obtain uniform bounds on $\vt$, 
\begin{equation} \label{kon}
0 < \underline{\vt} (T) \leq \vt(t,x) \leq \Ov{\vt}(T) \ \mbox{for all}\ t \in [0,T], \ x \in \Omega
\end{equation}
as soon as 
\[
0 < \inf_\Omega \vt_0 \leq \sup_\Omega \vt_0 < \infty.
\]

\subsubsection{Parabolic estimates}

Denoting 
\[
\kappa_\delta (\vt) = \kappa(s) + \delta \left( s^\Gamma + \frac{1}{s} \right),\ 
\mathcal{K}_\delta (\vt) = \int_1^\vt \left( \kappa(s) + \delta \left( s^\Gamma + \frac{1}{s} \right) \right) \D s
\]
we can rewrite \eqref{eE5} as
\begin{equation} \label{equ} 
\partial_t (\vr e_\delta) + \Div (\vr e_\delta \vu) - \Delta \mathcal{K}_\delta (\vt) =
\mathbb{S}_\delta  : \Grad \vu - p \Div \vu + \ep \delta \left( \Gamma \vr^{\Gamma - 2} + 2 \right) |\Grad \vr|^2 
+ \delta \frac{1}{\vt^2} - \ep \vt^5 
\end{equation}

Similarly to \cite[Chapter 3, Lemma 3.3]{FeNo6A}, the parabolic estimates 
\[
\sup_{t \in (0,T)} \| \vt \|_{W^{1,2}(\Omega)} + \int_0^T \intO{ \left( |\partial_t \vt|^2 + |\Del \mathcal{K}_\delta (\vt) |^2 
\right) } \dt \aleq 1
\]
can be derived by multiplying \eqref{equ} successively on $\vt$, $\partial_t \mathcal{K}_\delta (\vt)$. In the present setting, this 
technique produces two extra boundary integrals: 
\[
I_1 = \int_{\partial \Omega} \vt \left( \vr e_\delta \vuB \cdot \vc{n} + \vc{q}_\delta\cdot \vc{n} \right) \D \sigma_x = 
\int_{\Gamma_{\rm out}} \vt \vr e_\delta \vuB \cdot \vc{n} \D \sigma_x + \int_{\Gamma_{\rm in}} \vt F_{i,b} \ \D \sigma_x,
\] 
and 
\[
I_2 = \int_{\partial \Omega} \partial_t \mathcal{K}_\delta (\vt) \vc{q}_\delta \cdot \vc{n} \D \sigma_x = 
\frac{{\rm d}}{\dt} \int_{\Gamma_{\rm in}} \mathcal{K}_\delta (\vt) F_{i,b} \D \sigma_x - 
\int_{\Gamma_{\rm in}} \partial_t \mathcal{K}_\delta (\vt) \vr e_\delta (\vr, \vt) \vuB \cdot \vc{n} \D \sigma_x.
\]
Introducing a function  
\[
\chi (\vr, \vt) = \vr \int_1^\vt \kappa_\delta (s) e_\delta (\vr, s) \D s 
\]
we compute 
\[
\partial_t \chi (\vr, \vt) = \partial_t \vr \left( \int_1^\vt \kappa_\delta (s) e_\delta (\vr, s) \D s + 
\int_1^\vt \kappa_\delta (s) \frac{\partial e_\delta}{\partial \vr} (\vr, s) \D s \right) 
+ \partial_t \mathcal{K}(\vt) \vr e_\delta.  
\]
Consequently, 
\[
\begin{split}
I_2 &= \frac{{\rm d}}{\dt} \int_{\Gamma_{\rm in}} \Big( \mathcal{K}_\delta (\vt) F_{i,b} - \chi(\vr, \vt) \vuB \cdot \vc{n} \Big) \D \sigma_x
\\
&+ \int_{\Gamma_{\rm in}} \partial_t \vr \left( \int_1^\vt \kappa_\delta (s) e_\delta (\vr, s) \D s + 
\int_1^\vt \kappa_\delta (s) \frac{\partial e_\delta}{\partial \vr} (\vr, s) \D s \right) \D \sigma_x
\end{split}
\]
Thus all integrals are controlled in terms of $\partial_t \vr$ and the uniform bounds for $\vt$ established in \eqref{kon}.

\subsubsection{Solvability of the approximate internal energy equation}

Having established the same set of {\it a priori} bounds as in \cite[Chapter 3, Section 3.4.2]{FeNo6A}, the existence of the 
approximate solutions satisfying \eqref{eE5}, \eqref{eE6} can be shown as therein. Note that the Neumann problem for 
general quasilinear parabolic equations in divergence form is nowadays well understood. The relevant existence result was shown 
by Ladyzhenskaya, Solonnikov, and Uraltseva \cite[Chapter 5, par. 7, Theorem 7.4]{LSU} under certain restrictions imposed on the growth of the nonlinearities. As the comparison principle holds for the present problem, solutions may be constructed by suitable 
cut--off of nonlinearities, application of the result from \cite{LSU}, and passing to the limit in the regularization, see 
\cite[Chapter 3, Section 3.4.2]{FeNo6A} for details.

\subsection{Asymptotic limit of the approximate solutions}

Our ultimate goal is to show that any limit of a sequence of approximate solutions represents a weak solution of the target problem. 
Note that this includes limits at three different levels: 
\begin{itemize}
\item The limit $n \to \infty$ in the Galerkin approximation \eqref{E4} of the momentum equation. 
\item The artificial viscosity limit $\ep \to 0$ in the equation of continuity \eqref{E3}. 
\item The limit $\delta \to 0$ in the artificial regularizing terms.
\end{itemize}

A detailed proof of convergence is rather lengthy but nowadays well understood at least in the case of conservative boundary conditions. 
Indeed the full proof of convergence under the hypotheses \eqref{ws1}--\eqref{ws10} and with $\vuB \equiv 0$ was given in 
\cite[Chapter 3]{FeNo6A}. Moreover, the \emph{barotropic} Navier--Stokes system with general inflow/outflow boundary conditions has been treated in detail in \cite{ChJiNo}. Consequently, we focus only on the convergence of the boundary integrals in the total energy balance 
\eqref{E11}, and the entropy inequality \eqref{E8}. As the difficulties are the same at any level of the approximate process,  
we use a generic notation $(\vr_m, \vu_m, \vt_m)$ for an approximate sequence, where $m \to \infty$ stands for $n \to \infty$, or 
$\ep \to 0$, or $\delta \to 0$. We also focus on the last step $\delta \to 0$.

\subsubsection{Total energy balance} 

Neglecting all non--negative terms on the left--hand side of the approximate total energy balance {(\ref{E11})} we get 
\begin{equation}\label{E*}
\begin{split}
&\left[ \intO{\left[ \frac{1}{2} \vr |\vu - \vuB|^2 + \delta \left( \frac{1}{\Gamma - 1} \vr^\Gamma + \vr^2 \right) + 
\vr e_\delta  \right] \psi } \right]_{t = 0}^{ t = \tau}\\ &
- \int_0^\tau \partial_t \psi \intO{\left[ \frac{1}{2} \vr |\vu - \vuB|^2 + \delta \left( \frac{1}{\Gamma - 1} \vr^\Gamma + \vr^2 \right) + 
\vr e_\delta  \right] } \dt   
 \\ 
&+\int_0^\tau \psi \int_{\Gamma_{\rm out}} \vr e \vuB \cdot \vc{n}\  \D \sigma_x \dt    \dt 
\end{split}
\end{equation}
\begin{equation} \label{E17}
\begin{split}
&\leq - 
\int_0^\tau \psi \intO{ \left[ \vr \vu \otimes \vu + p_\delta  \mathbb{I} \right]  :  \Grad \vu_B } \dt + \frac{1}{2} \int_0^\tau 
\psi {\intO{ {\vr} \vu  \cdot\Grad |\vu_B|^2  } }
\dt  
\\ &+ \int_0^\tau \psi \intO{ \mathbb{S}_\delta(\vt,  \Grad \vu) : \Grad \vuB } \dt + 
\int_0^\tau \intO{ \vr \vc{g} \cdot (\vu - \vuB) } \dt - \int_0^\tau \psi \int_{\Gamma_{\rm in}} F_{i,b} \D \sigma_x \dt\\ &+\delta \int_0^\tau \psi \intO{  \frac{1}{\vt^2} } \dt - \delta \int_0^\tau \psi \int_{\Gamma_{\rm in}} \frac{1}{\Gamma - 1}\vr_b^\Gamma  \vuB \cdot \vc{n} \ \D \sigma_x \dt
\end{split}
\end{equation} 
for any $\psi \in C^1[0,T]$, $\psi \geq 0$.
Moreover, 
\[
\delta \int_0^\tau \int_{\Gamma_{\rm in}} \frac{1}{\Gamma - 1}\vr_b^\Gamma  \vuB \cdot \vc{n} \ \D \sigma_x \dt 
\to 0 \ \mbox{as}\ \delta \to 0, 
\]
and, by virtue of the uniform bounds \eqref{E16}, 
\[
\delta \int_0^\tau \intO{  \frac{1}{\vt^2} } \dt \to 0 \ \mbox{as}\ \delta \to 0.
\]

As the convergence of the volume integrals on the right-hand side of \eqref{E17} was established in 
\cite[Chapter 3]{FeNo6A}, it remains to handle the boundary term
\[
\int_0^\tau \psi \int_{\Gamma_{\rm out}} \vr e \vuB \cdot \vc{n}\  \D \sigma_x \dt    \dt
\]
Using the last estimate in \eqref{E15} we deduce 
\[
\int_0^T \int_{\Gamma_{\rm out}} \Big[ (\vr e - \vr \Ov{\vt} s)  \Big] |\vuB \cdot \vc{n}| \D \sigma_x \dt
\aleq 1.  
\]
As the function $(\vr e) = E_{\rm int}(\vr, S)$ is convex in the conservative entropy variables $(\vr, S=\vr s)$ 
and $e$ is given by \eqref{ws2}, we get 
\[
\int_0^T \int_{\Gamma_{\rm out}} \left( |\vr s| + \vr^{5/3} \right) |\vuB \cdot \vc{n}| \D \sigma_x \dt \aleq 1.
\]
Consequently, given an approximating sequence $\vr_m$, $S_m = \vr_m S_m$, we may suppose, extracting a suitable subsequence as the case may be, that
\[
\begin{split}
\vr_m &\to \vr \ \mbox{weakly in} \ L^{5/3}((0,T) \times \Gamma_{\rm out}; |\vuB \cdot \vc{n}| \dx ),\\ 
S_m &\to^b S\ \mbox{in the biting sense in}\ L^{1}((0,T) \times \Gamma_{\rm out}; |\vuB \cdot \vc{n}| \dx ),\\ 
(\vr_m, S_m) \ &\mbox{generates a Young measure}\ \nu_{t,x},\ (t,x) \in (0,T) \times \Gamma_{\rm out},\\
S(t,x)& = \left< \nu_{t,x}; \hat S \right>:=\int_{R^2}
\hat S{\rm d}\nu_{t,x}(\hat S),
\end{split}
\]
see Ball \cite{BA} and Ball, Murat \cite[Section 3]{BAMU}.
In order to pass to the limit in the total energy balance (\ref{E*}), we have to show
\bFormula{E**}
\liminf_{m \to \infty} \int_0^\tau \int_{\Gamma_{\rm out}} \vr_m e(\vr_m, S_m) \vuB \cdot \vc{n} \ \D \sigma_x \dt
\eF
$$
\geq \int_0^\tau \int_{\Gamma_{\rm out}} E_{\rm int} (\vr, S) \vuB \cdot \vc{n} \D \sigma_x \dt.
$$

As the function $E_{\rm int}$ being non--negative lower semi--continous, it can be approximated, by virtue of Baire's theorem, 
by a sequence of continuous, compactly supported functions, 
\[
E_n \in C_c(R^2), \ E_n \geq 0,\ E_n \nearrow E_{{\rm int}}\ \mbox{pointwise}. 
\]
This yields 
\[
\liminf_{{m \to \infty}} \int_0^\tau \int_{\Gamma_{\rm out}} E_{\rm {int}}(\vr_m, S_m ) \vuB \cdot \vc{n} \D \sigma_x \dt 
\geq \int_0^\tau \int_{\Gamma_{\rm out}} \left< \nu_{t,x}; E_{\rm {int}}(\hat \vr, \hat S ) \right>  \vuB \cdot \vc{n} \D \sigma_x \dt,  
\]
{in view of the fact that}
$$
\lim_{m \to \infty} \int_0^\tau \int_{\Gamma_{\rm out}} E_n(\vr_m, S_m ) \vuB \cdot \vc{n} \D \sigma_x \dt 
{=}\ \int_0^\tau \int_{\Gamma_{\rm out}} \left< \nu_{t,x}; E_n(\hat \vr, \hat S ) \right>  \vuB \cdot \vc{n} \D \sigma_x \dt.
$$

Finally, Jensen's inequality yields the desired conclusion, 
\[
\int_0^\tau \int_{\Gamma_{\rm out}} \left< \nu_{t,x}; E_{\rm {int}}( \hat \vr, \hat S ) \right>  \vuB \cdot \vc{n} \D \sigma_x \dt
\geq \int_0^\tau \int_{\Gamma_{\rm out}} E_{\rm {int}}(\vr, S )  \vuB \cdot \vc{n} \D \sigma_x \dt,
\]
}


\begin{Remark} \label{RR1}

To justify the above arguments we must extend carefully the function $E_{\rm int}(\vr, S)$ outside its natural domain 
$\vr > 0$. First, the entropy $s$ admits a limit 
\[
s(\vr, \vt) \to \underline{s} \in \{ 0, - \infty \} \ \mbox{as}\ \vt \to 0 \ \mbox{for any fixed}\ \vr > 0.
\]

The case $\underline{s} = 0$ corresponds to the Third law of thermodynamics  and has been 
considered as one of the hypotheses in Theorem \ref{RT2}. In this case, the internal energy is defined as
\[
E_{\rm int}(\vr, S) = \left\{ \begin{array}{l} \vr e(\vr, S) \ \mbox{if}\ \vr > 0, \ S > 0,\\ \\ 
0 \ \mbox{if}\ \vr = 0,\ S = 0, \\ \\ \lim_{\vr \to 0+} \vr e(\vr, S) \ \mbox{if}\ \vr = 0, \ S > 0,\\ \\ 
\infty \ \mbox{otherwise.} \end{array} \right.
\] 

If $\underline{s} = - \infty$ we set 
\[
E_{\rm int}(\vr, S) = \left\{ \begin{array}{l} \vr e(\vr, S) \ \mbox{if}\ \vr > 0, \ S \in R,\\ \\ \lim_{\vr \to 0+} \vr e(\vr, S) \ \mbox{if}\ \vr = 0, \ S \in R,\\ \\ 
\infty \ \mbox{otherwise.} \end{array} \right.
\] 

In both cases the extended function is convex, l.s.c, and strictly convex in the interior of its domain.

\end{Remark}

\subsubsection{Entropy inequality}

To finish the proof of convergence, we have to perform the limit in the boundary integrals in the approximate entropy inequality 
\eqref{E8}. Let us start with 
\begin{equation}\label{E***}
\int_0^\tau \int_{\Gamma_{\rm out}} \varphi \vr s_\delta \vuB \cdot \vc{n} \ \D \sigma_x \dt = \delta \int_{\Gamma_{\rm out}} \varphi \vr \log(\vt) \vuB \cdot \vc{n} \ \D \sigma_x \dt +
\int_0^\tau \int_{\Gamma_{\rm out}} \varphi \vr s \vuB \cdot \vc{n} \ \D \sigma_x \dt.
\end{equation}

{We start with the first term in (\ref{E***}). We have,
\[
\delta \int_{\Gamma_{\rm out}} \varphi \vr \log(\vt) \vuB \cdot \vc{n} \ \D \sigma_x \dt \leq 
\delta \int_{\Gamma_{\rm out}} \varphi \vr [\log \vt]^+ \vuB \cdot \vc{n} \ \D \sigma_x \dt \to 0 
\ \mbox{as}\ \delta \to 0
\]
as a consequence of uniform integrability of the internal energy. 

The treatement of the second term in (\ref{E***}) is more delicate.
We have to show that 
\begin{equation} \label{E18}
\limsup_{m \to \infty} \int_0^\tau \int_{\Gamma_{\rm out}} \varphi S_m \vuB \cdot \vc{n} \ \D \sigma_x \dt 
\leq \int_0^\tau \int_{\Gamma_{\rm out}} \varphi S \vuB \cdot \vc{n} \ {\D \sigma_x \dt, \
\varphi \geq 0,}
\end{equation}
where $S$ is the biting limit of the sequence $\{ S_m \}_{m=1}^\infty$.
To this end consider a function 
\[
\chi \in C^\infty(R), \ 0 \leq \chi \leq 1, \ \chi' \geq 0 , \ \chi(Y) = 0 \ \mbox{for}\ Y \leq 0,\ 
\chi(Y) = 1 \ \mbox{for}\ Y \geq 1,
\]
and the composition 
\[
\chi(S + k) ,\ k > 1.
\]
Observe that
\[
S \leq \chi(S + k) S \ \mbox{if}\ k > 1.
\]
Consequently, 
\[
\limsup_{m \to \infty} \int_0^\tau \int_{\Gamma_{\rm out}} \varphi S_m \vuB \cdot \vc{n} \ \D \sigma_x \dt 
\leq \limsup_{m \to \infty} \int_0^\tau \int_{\Gamma_{\rm out}} \varphi \chi (S_m + k) S_m \vuB \cdot \vc{n} \ \D \sigma_x \dt
\ \mbox{for any}\ k > 1.
\]
Next, we claim that the family $\{ \chi (S_m + k) S_m \}_{m > 0}$ is equi--integrable in 
$L^1((0,T) \times \Gamma_{\rm out}; |\vuB \cdot \vc{n}|\dx )$ for any fixed $k > 1$. Assuming for a moment this is the case, 
we get 
\[
\begin{split} 
&\limsup_{m \to \infty} \int_0^\tau \int_{\Gamma_{\rm out}} \varphi S_m \vuB \cdot \vc{n} \ \D \sigma_x \dt 
\leq \limsup_{m \to \infty} \int_0^\tau \int_{\Gamma_{\rm out}} \varphi \chi (S_m + k) S_m \vuB \cdot \vc{n} \ \D \sigma_x \dt\\
&= \int_0^\tau \int_{\Gamma_{\rm out}} \varphi \left< \nu_{t,x}; \chi ({\hat S} + k) \hat S \right> \vuB \cdot \vc{n} \ \D \sigma_x \dt
\to \int_0^\tau \int_{\Gamma_{\rm out}} \varphi \left< \nu_{t,x}; \hat S \right> \vuB \cdot \vc{n} \ \D \sigma_x \dt
\ \mbox{as}\ k \to \infty,
\end{split}
\]
as claimed in \eqref{E18}. To see equi--integrability of $\{ \chi (S_m + k) S_m \}_{m > 0}$, consider first 
the part of $(0,T) \times \Gamma_{\rm out}$, where $\vt_m > 1$. By virtue of the hypotheses \eqref{ws4}--\eqref{ws6}, 
\[
|S_m| = |\vr_m s(\vr_m, \vt_m)| \aleq \left( 1 + \vr_m |\log(\vr_m)| + (\vt^m)^3 \right) 
\ \mbox{if}\ \vt_m > 1 
\]
and the desired equi--integrability follows from uniform integrability of $E_{\rm in}(\vr_m, \vt_m)$. If $\vt_m \leq 1$, 
we deduce from \eqref{ws3} that 
\[
S_m \aleq \vr_m \mathcal{S}(\vr_m) + 1 \aleq E_{\rm in}(\vr_m, \vt_m);
\]
whence 
\[
- (k+1) \leq \chi (S_m + k) S_m \leq E_{\rm in}(\vr_m, \vt_m)
\]
which implies equi--integrability. {This finishes the proof of (\ref{E18}).}


Thus it remains to pass to the limit in the boundary integral
\[
\int_0^\tau \int_{\Gamma_{\rm in}} \varphi \left[ \delta \vr_b (1 - \log(\vt) ) 
\vuB \cdot \vc{n} - \frac{1}{\vt} F_{i,b} + 
\left( \frac{e(\vr_b, \vt)}{\vt} - s(\vr_b, \vt) \right) \vr_b \vuB \cdot \vc{n}  \right] \D \sigma_x \dt. 
\]
By virtue of the uniform bounds \eqref{E15}, we have 
\[
\int_0^T \intO{ \| \log(\vt) \|^2_{W^{1,2}(\Omega)} } \dt \aleq 1; 
\]
whence, as a consequence of the trace theorem, 
\[
\delta \int_0^\tau \int_{\Gamma_{\rm in}} \varphi \vr_b (1 - \log(\vt) ) 
\vuB \cdot \vc{n} \ \D \sigma_x \dt \to 0 \ \mbox{as}\ \delta \to 0.
\]
Next, exploiting \eqref{E15} once more we get 
\[
{\rm ess} \sup_{t \in (0,T)} \| \vt \|_{L^4(\Omega)} + \int_0^T \|\Grad \vt \|^2_{L^2(\Omega)} \aleq 1.
\]
Moreover, 
as shown in \cite[Chapter 3]{FeNo6}, 
\[
\vt_m \to \vt \ \mbox{in, say,}\ L^2((0,T) \times \Omega).
\]
Thus, by interpolation, 
\[
\int_0^T \int_{\partial \Omega} \| \vt_m - \vt \|^2 \ \D \sigma_x \dt \aleq \int_0^T 
\| \vt_m - \vt \|^2_{W^{\alpha, 2}(\Omega) } \dt \leq 
\int_0^T 
\| \vt_m - \vt \|^{2\alpha}_{W^{1, 2}(\Omega) }\| \vt_m - \vt \|_{L^2(\Omega) }^{2 (1 - \alpha)}   \dt
\]
for any $\frac{1}{2} < \alpha \leq 1$. Consequently, 
\[
\vt_m \to \vt \ \mbox{in}\ L^2((0,T) \times \partial \Omega), 
\]
and the limit in the integral 
\[
\int_0^\tau \int_{\Gamma_{\rm in}} \varphi \left[ - \frac{1}{\vt_m} F_{i,b} + 
\left( \frac{e(\vr_b, \vt_m)}{\vt_m} - s(\vr_b, \vt_m) \right) \vr_b \vuB \cdot \vc{n}  \right] \D \sigma_x \dt, 
\]
can be performed using hypothesis \eqref{E1} and Fatou's lemma.

We have shown Theorem \ref{tE1}.

\section{Concluding remarks} 

To the best of our knowledge, this is the first result concerning global existence for the Navier--Stokes--Fourier system with 
large and realistic initial/boundary conditions. The fact that the present concept of weak solution complies with the weak--strong uniqueness 
principle plays an important role, in particular in view of the recent results on ill--posedness of the incompressible Navier--Stokes system, see Buckmaster and Vicol \cite{BucVic}. The present results can be used as a suitable platform for studying turbulence phenomena in physically relevant open fluid systems.  

Extensions to more general rheological laws are certainly possible, however, the basic structure of the internal energy
\[
\vr e(\vr, \vt) \approx \vr \vt + \vr^\gamma + \vt^4 
\] 
is essential in view of the lack of suitable {\it a priori} bounds. In particular, the stabilizing effect of the radiation 
component $\vt^4$ is absolutely crucial on (hypothetical) vacuum zones, where $\vr$ vanishes.

Last but not the least, the approximate scheme used in the construction of weak solutions shares certain similarity with the numerical methods based on the upwinding of convective terms, see 
\cite{FeKaNo},  \cite{FeiLMHMShe}, \cite{KwNo}. 

\def\cprime{$'$} \def\ocirc#1{\ifmmode\setbox0=\hbox{$#1$}\dimen0=\ht0
  \advance\dimen0 by1pt\rlap{\hbox to\wd0{\hss\raise\dimen0
  \hbox{\hskip.2em$\scriptscriptstyle\circ$}\hss}}#1\else {\accent"17 #1}\fi}

\end{document}